\documentclass[preprint,showpacs,english,aps]{revtex4-1}

\usepackage{babel,rotating,dcolumn}
\usepackage{colordvi,graphicx,color,amsbsy,amsmath,bm,amssymb}
\usepackage{comment}

\usepackage{graphicx}
\usepackage{epstopdf, epsfig}

\usepackage{babel,rotating,dcolumn}
\usepackage{colordvi,graphicx,color,amsbsy,amsmath,bm,amssymb}
\usepackage{comment}
\usepackage[normalem]{ulem}

\DeclareSymbolFont{matha}{OML}{txmi}{m}{b}
\DeclareMathSymbol{\varv}{\mathord}{matha}{118}
\usepackage{mathrsfs}

\DeclareMathAlphabet\mathbfcal{OMS}{cmsy}{b}{n}

\usepackage{yhmath}
\usepackage{ucs}
\usepackage{amsfonts}

\usepackage{amsmath,esint}

\usepackage{nomencl}

\usepackage{bbm}

\usepackage[dvipsnames]{xcolor}

\newcommand*\Thibault[1]{\textcolor{ForestGreen}{#1}}
\newcommand*\Keigo[1]{\textcolor{red}{#1}}
\newcommand*\Kai[1]{\textcolor{blue}{#1}}

\begin{document}

\title[]{Computing differential operators of the particle velocity in moving particle clouds using tessellations}
\author{Thibault Maurel--Oujia$^{1}$}
\email{thibault.oujia@univ-amu.fr}
\author{Keigo Matsuda$^{2}$} 
\author{Kai Schneider$^{1}$}
%
\affiliation{$^{1}$ Institut de Math\'ematiques de Marseille, Aix-Marseille Université, CNRS, France}
\affiliation{$^{2}$ Japan Agency for Marine-Earth Science and Technology (JAMSTEC), Japan}
\date{\today}

\begin{abstract}
We propose finite-time measures to compute the divergence, the curl and the velocity gradient tensor of the point particle velocity for two- and three-dimensional moving particle clouds. 
For this purpose, a tessellation of the particle positions is performed to assign a volume to each particle.
We introduce a modified Voronoi tessellation which overcomes some drawbacks of the classical construction.
Instead of the circumcenter we use the center of gravity of the Delaunay cell for defining the vertices.  
Considering then two subsequent time instants, the dynamics of the volume can be assessed. 
Determining the volume change of tessellation cells yields the divergence of the particle velocity.
Reorganizing the various velocity coefficients allows computing the curl and even the velocity gradient tensor.
The helicity of particle velocity can be likewise computed and swirling motion of particle clouds can be quantified. 
First we assess the numerical accuracy for randomly distributed particles. 
We find a strong Pearson correlation between the divergence computed with the the modified tessellation, and the exact value.
Moreover, we show that the proposed method converges with first order in space and time in two and three dimensions.
Then we consider particles advected with random velocity fields with imposed power-law energy spectra.
We study the number of particles necessary to guarantee a given precision. 
Finally, applications to fluid particles advected in three-dimensional fully developed isotropic turbulence show the utility of the approach for real world applications to quantify self-organization in particle clouds and their vortical or even swirling motion.
\end{abstract}

\maketitle

\section{Introduction}
The general context of this work are flows of point particle clouds evolving in time where inertial particles are transported by an underlying flow velocity, the so-called the carrier phase. Such configurations are typically encountered in multiphase flows, especially when the carrier phase velocity is turbulent. 
Particle laden turbulent flows occur in many natural and engineering systems, for instance droplets in clouds in the atmosphere, sand storms in the desert or spray combustion with liquid fuels. 
Observations show that inertial particles in turbulent flow exhibit a non-uniform spatial distribution, i.e., preferential concentration, resulting in the emergence of void and cluster regions \citep{squires1990particle, squires1991preferential}. %
Understanding of the complex dynamics of these fundamental cluster and void formation processes is a prerequisite for sound modeling and important for numerous applications. 
The divergence of the particle velocity, a scalar-valued signed quantity, allows to quantify the divergence and convergence of particles corresponding to positive and negative divergence values, respectively. 
The divergence of the particle velocity is contained in the source term in the particle density balance equation. 
The curl of the particle velocity, a vector valued quantity, yields insight about the local rotation rate of particles, i.e., the particle vorticity and thus quantifies vortical motion of particle clouds. Likewise important to characterize swirling motion is the flow helicity which quantifies the alignment of the particle velocity and the particle vorticity.
The ultimate goal is to access the particle velocity gradient tensor, and to consider its trace, the symmetric and anti-symmetric parts. 
%
%
%

However, it is not straightforward
to determine quantities such as divergence and curl of the particle velocity due to its discrete nature, i.e., the velocity is only defined at particle positions. 
One approach to deal with the problem for instance in homogeneous periodic flow would be to
determine the particle velocity from the particle flux.
The particle flux is given by the convolution of the particle density with the particle velocity, and in principle we could compute the particle velocity by solving a deconvolution problem.
However, the particle positions are not given on a uniform grid. Hence the Fast Fourier Transform (FFT) cannot be used, which makes this problem very expensive and even unfeasible for a large number of particles.
Interpolating the particle velocity data on a regular grid, using e.g. Fourier, or other interpolation methods would be another option. Then the derivative of the particle velocity could be computed by deriving the interpolating function. However, interpolating values given on random particle positions has its own difficulties.
%
%
For solving PDEs and also for computing differential operators different meshfree methods have been developed, some review detailing the major directions, i.e. smooth particle hydrodynamics and moving least-squares, is given in \citet{huerta2004meshfree}. 

Voronoi and Delaunay tessellations allow the definition of natural neighbors for arbitrary unstructured grids.
Voronoi tessellation has been applied to analyze data from different origins, e.g. astrophysics, biology, particle-laden turbulence \citep{ebeling1993detecting, OTCMB14, DeMo13}. 
The volume of Voronoi cells yields the local density of a set of points. 
Using Voronoi tessellation, statistics are not affected by an arbitrary choice, like the bin size, for defining for instance the particle density.
A mimetic method to compute divergence and curl operators on unstructured grids using Voronoi and Delaunay tessellations has been developed in \citet{vabishchevich2005finite}, by refining a velocity gradient on the Voronoi and Delaunay cells. This approach allows to define the different operators in two and three space dimensions. 
However, using this method the truncation error was shown to be in general a zeroth-order accuracy, except for regular grids. 
\citet{sozer2014gradient} assessed the accuracy of gradient reconstruction methods for cell-centered data on unstructured meshes. 
A Green--Gauss method is used for gradient calculation with different interpolation stencils (e.g.  Weighted Least Squares).
It is shown that several gradient operators have first order accuracy, regardless of the cell shape. 
However, there is not one gradient operator which yields 
the lowest error for all cell configurations. 
An analysis of the truncation error on a structured and unstructured grid of the gradient calculus using the Green--Gauss method and least-squares method can be found in \citet{syrakos2017critical}.
%
%
%
%

%
%
%
In \citet{oujia_matsuda_schneider_2020} we proposed a finite--time measure to quantify the divergence, by combining tessellation with a Lagrangian approach. The divergence was computed by determining the volume change rate of the Voronoi volume and we presented some preliminary validation. 

The aim in the present work is to propose efficient numerical methods to compute differential operators of the particle velocity using volume change of the tessellation cells to asses the dynamics of moving particle clouds.
In this article we introduce an improved method and compare it to the one presented in \citet{oujia_matsuda_schneider_2020} for calculating the divergence of a particle velocity field. 
Moreover, we extend this approach for the computation of the curl and velocity gradient tensor of the particle velocity. 
In addition, we perform a detailed numerical analysis of the properties and convergence of the different methods. 
Furthermore, we determine the minimum number of particles required to achieve a given level of accuracy in a turbulent flow in two and three dimensions. The results of this study are of practical importance for researchers interested in accurately computing velocity derivatives in particle-based simulations and to quantify the dynamics of clustering.
The outline of the manuscript is as follows. In section~\ref{sec:mathematical_foundation} we summarize the principle of Voronoi and Delaunay tessellation and explain the mathematical foundation of the proposed approach to quantify the divergence, curl and velocity gradient tensor of particle velocity using volume change. 
In section~\ref{sec:numerical_validation} we provide some information on the reliability and analyze the convergence order of our methods. 
In section~\ref{sec:application_to_turbulence} we apply the methods to fluid and inertial particles driven by turbulent velocity fields. Finally, conclusions are drawn in section~\ref{sec:conclusions}. 
\section{Differential operators based on volume change}
\label{sec:mathematical_foundation}

In this section we explain the mathematical foundation of our method to define 
finite–time measures to quantify spatial differential operators.
First, Voronoi and Delaunay tessellations of particle clouds are described, including some mathematical background.
Then, we present how to calculate the divergence from the volume change over time and adapt the method to enable 
the calculation of the curl and the velocity gradient tensor of the particle velocity.
Finally, we convert this Lagrangian method into an Eulerian method to estimate the truncation error. 

\subsection{Modified Voronoi tessellations}

The Voronoi tessellation~\citep{voronoi1908nouvelles} decomposes the space seeded with particles into cells; a Voronoi cell $\mathscr{C}_i$  
generated by the particle $p_i$ represents the region of space closer to the particle $p_i$ than to all the others.
A Voronoi cell can be interpreted as the influence volume of a particle.
The volumes of the Voronoi cells are inversely proportional to the particle density. 
For instance, we can use the size of the cells to identify clusters of particles. 
Indeed, a set of small cells is the sign of a local concentration of particles and a set of large cells indicates a low density region.
A quality of the Voronoi tessellation is that the tessellation is constructed without introducing an arbitrary length scale.
The Delaunay tessellation~\citep{delaunay1934sphere} of a set of particles $p$ is a tessellation such that no particles are inside the circumscribed circle of one of any cells. 
The Delaunay tessellation of a set of particles $p$ is the dual graph of the Voronoi tessellation. The vertices of the Voronoi tessellation are the centers of the circumscribed circles of the cells of the Delaunay tessellation. 
Another way to define a cell corresponding to a particle, instead of using the circumcenter of the Delaunay cell as done for the Voronoi tessellation, is to use the center of gravity of the Delaunay cell to define the vertices of the cell. 
In the following, the tessellation using the center of gravity of the Delaunay cell is called modified Voronoi tessellation.
We will see later that this construction yields more accurate results for computing differential operators.
To compute the Voronoi and the Delaunay tessellation, we use the Quickhull algorithm provided by the Qhull library in python \citep{BaDH96}, which has a computational complexity of ${\cal O}(N \log(N))$, where $N$ denotes the number of particles. 

\begin{figure}
\centering
\includegraphics[width=0.90\linewidth]{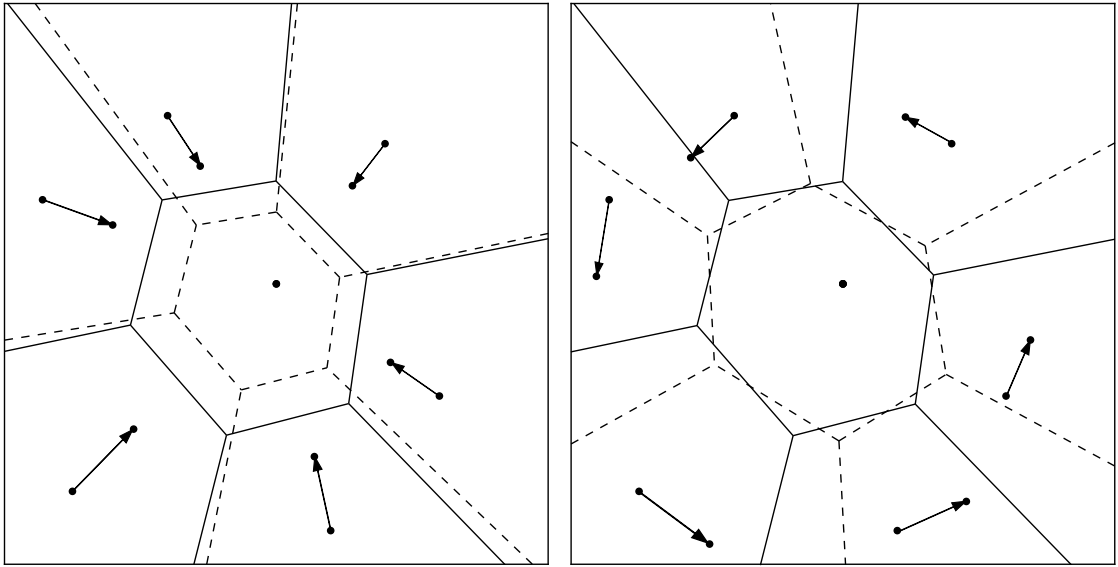}
\caption{Example of the motion of a two-dimensional Voronoi tessellation at two subsequent time instants for a compression (left) and a rotation (right). The solid line corresponds to the time step $t^k$ and the dashed line to $t^{k+1}$. The arrows indicate the particle velocity. 
} 
\label{fig:Illust_Div_Curl}
\end{figure}
Figure~\ref{fig:Illust_Div_Curl} illustrates a possible evolution of a two-dimensional Voronoi cell over time. On the left figure we can observe a compression of the Voronoi cell, while on the right we can see a rotation of the cell.
We can discern that in the case where the particles converge to the center, the volume of the Voronoi cell is reduced, and that in the case of rotation the cell is likewise rotated.
This figure illustrates the underlying idea of the method, which is to observe the time the evolution of volumes.
In this study, we focus on two and three dimensions, but the proposed methods can be directly applied to higher dimensions.

\subsection{Divergence of particle velocity using volume change}
\label{sec:div}

In this subsection, we first recall the method proposed in \citet{oujia_matsuda_schneider_2020} 
to define the finite time--discrete divergence of the particle velocity based on the time change of the Voronoi volume.
Motivated by stability concerns this method will be here also extended using modified Voronoi volumes. 

The particles are transported by a velocity field. Thus we have access to the particle velocity at the discrete particle position.
Particles in a 
flow domain $\Omega \subset \mathbb{R}^d$, $d=2,3$, satisfy the conservation equation, 
\begin{equation} \label{eq:conservation}
D_t n = \partial_t n + {\bm v}\cdot\nabla n = - n \nabla \cdot {\bm v}, 
\end{equation}
as already shown in
\citet{Maxe87}, 
where $n \ge 0$ is the particle number density, ${\bm v}$ is the particle velocity, and $D_t n$ is the Lagrangian derivative of $n$. 
We use the notation $\partial_t=\partial/\partial t$ and $\nabla = (\partial_x, \partial_y, \partial_z)^T$ where the superscript $^T$ is the transposition operator.
%
Dividing by the particle density $n$ and considering the particle density at two time instants $t^k$ and $t^{k+1}=t^k+\Delta t$,
\citet{oujia_matsuda_schneider_2020} define the 
divergence 
$\nabla \cdot {\bm v}$ using the particle density
as, 
\begin{equation}\label{eq:div_density}
\nabla \cdot {\bm v} =
-\frac{1}{n}D_t n = -\frac{2}{n^{k+1}+n^{k}} \, \frac{n^{k+1}-n^{k}}{\Delta t} + {\cal O}(\Delta t).
\end{equation}
%
%
%
%
To obtain the particle distribution at time $t^{k+1}$ particles are advected by the particle velocity, e.g., using the explicit Euler method, ${\bm x}^{k+1}_p = {\bm x}^k_p+\Delta t ~ {\bm v}_p^k$ where ${\bm x}_p$ and ${\bm v}_p$ are the position and velocity of the particles, respectively.

The local number density over a Voronoi cell $n_p>0$, a piecewise constant function, is defined as $n_p = V_{p}^{-1}$, with $V_{p}$ being the volume of the cell. 
Similarly, the local number density can be defined also based on modified Voronoi cells, i.e., defined based on centers of gravity instead of circumcenters. 
Replacing the particle density $n$ by this discrete value in equation (\ref{eq:div_density}), we obtain a finite-time discrete divergence of the particle velocity based on tessellation cells defined as
\begin{equation} \label{eq:div_density_tessellation}
{\cal D}({\bm v}_p) = \frac{2}{\Delta t} \, \frac{V_p^{k+1}-V_p^{k}}{V_p^{k+1}+V_p^{k}}
\end{equation}
where the subscript $_p$ represents the value of the function at the particle position ${\bm x}_p$. 
%

Other temporal approximations of the particle density $n$ are possible, e.g.  using linear time interpolation of
$1/n$, 
 which gives
\begin{equation}  \label{eq:div_density_tessellation_lin}
{\cal D}_{lin} ({\bm v}_p) = \frac{1}{2\Delta t}\left( \frac{1}{V_p^{k+1}} + \frac{1}{V_p^{k}} \right) \, (V_p^{k+1}-V_p^{k})
\end{equation}
or using the logarithmic derivative
\begin{equation}  \label{eq:div_density_tessellation_log}
{\cal D}_{log} ({\bm v}_p) 
= \frac{1}{\Delta t}\ln{\left(\frac{V_p^{k+1}}{V_p^{k}}\right)}.
\end{equation}
These three different approximations, which all have first order accuracy in time, are also tested in subsection \ref{sec:Computation_of_a_spatial_derivative_2D} in order to verify that the formula~(\ref{eq:div_density_tessellation}), presented in \citet{oujia_matsuda_schneider_2020}, gives the most accurate results for large values of $\Delta t$ and same results for small $\Delta t$ values.
All results shown in the following are thus based on equation (\ref{eq:div_density_tessellation}), unless otherwise specified.




\subsection{Rotation of particle velocity using volume change}
\label{sec:curl} 

The curl of the particle velocity 
can be defined by computing the circulation of the velocity field of particles over a cell. 
%
%
The curl in the direction of the $x$-axis, which is the first component of the curl, can be expressed as the negative 
divergence of a function ${\bm v}^\perp_x$, where ${\bm v}^\perp_x$ is obtained by projecting ${\bm v}$ orthogonally onto the plane perpendicular to the $x$-direction, 
followed by a counter-clockwise rotation of $\pi/2$ about the $x$-axis. The same procedure can be applied for the curl in the direction of the $y$- and $z$-axis,
%
%
%
%
\begin{equation}
    \nabla\times {\bm v} = 
    \begin{pmatrix}
    \displaystyle \partial_y v_z -  \partial_z v_y \\
    \displaystyle \partial_z v_x -  \partial_x v_z \\
    \displaystyle \partial_x v_y -  \partial_y v_x
    \end{pmatrix} = 
    \begin{pmatrix}
    0&-\displaystyle \partial_z &~ \partial_y  \\ 
    ~\displaystyle \partial_z &0&-\displaystyle \partial_x  \\ 
    -\displaystyle \partial_y &~\displaystyle \partial_x &0
    \end{pmatrix} 
    \begin{pmatrix}v_x \\ v_y \\ v_z\end{pmatrix}= 
    -\begin{pmatrix}\nabla\cdot{\bm v}^\perp_x\\ \nabla\cdot{\bm v}^\perp_y\\ \nabla\cdot{\bm v}^\perp_z\end{pmatrix},
    \label{eq:curl_3D_div_perp}
\end{equation}
where ${\bm v}=(v_x, v_y, v_z)^T$, 
%
$
{\bm v}^\perp_x = 
{\bm L}_x{\bm v}
$
, 
$
{\bm v}^\perp_y = 
{\bm L}_y{\bm v}
$
 and 
$
{\bm v}^\perp_z = 
{\bm L}_z{\bm v}
$ 
with \\
\begin{equation}
{\bm L}_x = 
\begin{pmatrix}
~0&~0&~0 \\ 
~0&~0&-1 \\ 
~0&~1&~0
\end{pmatrix}
\text{, }
{\bm L}_y = 
\begin{pmatrix}
~0&~0&~1 \\ 
~0&~0&~0 \\ 
-1&~0&~0
\end{pmatrix}
\text{ and }
{\bm L}_z = 
\begin{pmatrix}
~0&-1&~0 \\ 
~1&~0&~0 \\ 
~0&~0&~0
\end{pmatrix}.
\end{equation}
Note that $\{{\bm L}_x, {\bm L}_y, {\bm L}_z\}$ is the most common basis of the Lie algebra of the rotation group $SO(3)$ \citep{kirillov2008introduction}.
%
%
%
%
%
%
By reusing the discrete divergence ${\cal D}$ defined in subsection~\ref{sec:div}, the discrete curl $\mathbfcal{C}$ of the particle velocity can be computed at particle positions.
%
%
Using equation (\ref{eq:curl_3D_div_perp}) the discrete curl can be defined as
\begin{equation}
    \mathbfcal{C}({\bm v}_p)
    = \begin{pmatrix} {\cal D}\left(-{\bm v}^\perp_{p,x}\right)\\ {\cal D}\left(-{\bm v}^\perp_{p,y}\right)\\ {\cal D}\left(-{\bm v}^\perp_{p,z}\right) \end{pmatrix}. 
    \label{eq:curl_3D_div_perp_tessellation}
\end{equation}
In two dimensions, the curl reduces to a pseudo-scalar ${\cal C}({\bm v}_p)  = {\cal D}(-{\bm v}_p^\perp)$ where ${\bm v}_p^\perp$ is the particle velocity ${\bm v}_p$, rotated $\pi/2$ in counter-clockwise direction in the two-dimensional plane. 
Note that, the curl of the fluid velocity is the vorticity ${\bm \omega}$ and similarly the curl of the particle velocity can be considered to be the particle vorticity.
%

\subsection{Velocity gradient tensor of the particle velocity using volume change}
\label{sec:tensor} 

The velocity gradient tensor, 
can be computed as the divergence 
of each of its components projected onto the different axes of the vector space. 
In three dimensions, the velocity gradient tensor is defined and can be rewritten as,

\begin{equation}
\nabla{\bm v} = \begin{pmatrix}
  \displaystyle{\partial_x v_x} & \displaystyle{\partial_y v_x} & \displaystyle{\partial_z v_x} \\
  \displaystyle{\partial_x v_y} & \displaystyle{\partial_y v_y} & \displaystyle{\partial_z v_y} \\
  \displaystyle{\partial_x v_z} & \displaystyle{\partial_y v_z} & \displaystyle{\partial_z v_z}
\end{pmatrix}
=
\begin{pmatrix}
\nabla \cdot (v_x{\bm e}_1) & \nabla \cdot (v_x{\bm e}_2) & \nabla \cdot (v_x{\bm e}_3)\\
\nabla \cdot (v_y{\bm e}_1) & \nabla \cdot (v_y{\bm e}_2) & \nabla \cdot (v_y{\bm e}_3)\\
\nabla \cdot (v_z{\bm e}_1) & \nabla \cdot (v_z{\bm e}_2) & \nabla \cdot (v_z{\bm e}_3)
\end{pmatrix},
\label{eq:grad_v_3D}
\end{equation}
with
$$
{\bm e}_1 = \begin{pmatrix}1 \\0 \\0\end{pmatrix}; \;
{\bm e}_2 = \begin{pmatrix}0 \\1 \\0\end{pmatrix}; \;
{\bm e}_3 = \begin{pmatrix}0 \\0 \\1\end{pmatrix}
.
$$
%
Using the volume change technique to compute the discrete divergence of 
components of the function ${\bm v}_p$ projected on the different space axes, 
we can define the discrete gradient $\mathbfcal{G}$ of the velocity ${\bm v}_p$ 
at particle positions as
\begin{equation}
    \mathbfcal{G}({\bm v}_p) = 
    \begin{pmatrix}
    {\cal D}(v_{p,x}{\bm e}_1) & {\cal D}(v_{p,x}{\bm e}_2) & {\cal D}(v_{p,x}{\bm e}_3)\\
    {\cal D}(v_{p,y}{\bm e}_1) & {\cal D}(v_{p,y}{\bm e}_2) & {\cal D}(v_{p,y}{\bm e}_3)\\
    {\cal D}(v_{p,z}{\bm e}_1) & {\cal D}(v_{p,z}{\bm e}_2) & {\cal D}(v_{p,z}{\bm e}_3)
    \end{pmatrix}
    \label{eq:grad_v_3D_tessellation}.
\end{equation}
In a similar way, we can compute the velocity gradient tensor in two dimensions.

\subsection{Conversion into an Eulerian formulation and truncation error} 
\label{sec:euler}

The above finite-time measure of the Lagrangian derivative should be equivalent with an Eulerian method in the limit of $\Delta t \to 0$. Therefore, here we derive the Eulerian formula corresponding to the proposed Lagrangian approach and deduce the spatial accuracy as a truncation of the Taylor expansion.
Here, we focus on the three-dimensional case, similar results can be derived in two dimensions.
We consider four points $p_0$, $p_1$, $p_2$, and $p_3$. 
Here $p_0$ is considered to be the center particle, and $p_1$, $p_2$, and $p_3$ are considered to be vertices of a tessellation cell. 
The position and velocity vectors of $p_j$ ($j=0$, 1, 2, and 3) are ${\bm x}_j$ and ${\bm v}_j$, respectively.
We define the relative position vector ${\bm r}_j={\bm x}_j-{\bm x}_0$ and the relative velocity vector ${\bm w}_j={\bm v}_j-{\bm v}_0$ for $j=1$, 2, and 3.
The volume $V$ of the parallelepiped generated by the points $p_0$, $p_1$, $p_2$, and $p_3$ is given by $V=({\bm r}_1\times{\bm r}_2)\cdot{\bm r}_3$.
Note that the volume of the tetrahedron is 1/6 of the one of the parallelepiped.
Since the volume $V$ at the time $t+\Delta t$ is given by $V(t+\Delta t)=\{({\bm r}_1+{\bm w}_1\Delta t)\times({\bm r}_2+{\bm w}_2\Delta t)\}\cdot({\bm r}_3+{\bm w}_3\Delta t)$, the volume change of the parallelepiped is then given by 
\begin{equation}
\frac{d V}{d t}
=\lim_{\Delta t \to 0} \frac{V(t+\Delta t)-V(t)}{\Delta t} \\
=({\bm r}_2\times{\bm r}_3)\cdot{\bm w}_1 + ({\bm r}_3\times{\bm r}_1)\cdot{\bm w}_2 + ({\bm r}_1\times{\bm r}_2)\cdot{\bm w}_3
\end{equation}
Thus, assuming that the terms of ${\cal O}(\Delta t)$ are negligibly small, we 
obtain the divergence at the limit of $\Delta t \to 0$:
\begin{equation}
    {\cal D}_0({\bm v}) \equiv
    \frac{1}{V}\frac{d V}{d t} 
    = \frac{1}{V}[{\bm w}_1\cdot({\bm r}_2\times{\bm r}_3) 
     +{\bm w}_2\cdot({\bm r}_3\times{\bm r}_1) 
     +{\bm w}_3\cdot({\bm r}_1\times{\bm r}_2) ]
     \label{eq:div_euler}.
\end{equation}
Note that in 2D, ${\bm r}_3$ is replaced by the unit normal vector ${\bm n}={\bm r}_1\times{\bm r}_2/|{\bm r}_1\times{\bm r}_2|$, and ${\bm w}_3$ is replaced by ${\bm 0}$.

We consider the coordinate transformation in Cartesian coordinates of $(x_1,x_2,x_3)$ into $(\xi_1,\xi_2,\xi_3)$, in which a vector ${\bm x} = x_1{\bm e}_1 + x_2{\bm e}_2 + x_3{\bm e}_3$ is transformed into ${\bm x} = \xi_1{\bm r}_1 + \xi_2{\bm r}_2 + \xi_3{\bm r}_3$. 
The divergence of ${\bm v}$ in the transformed coordinates is then given by
\begin{equation}
    \nabla\cdot{\bm v} = {\bm e}^j\partial_{x_j} \cdot {\bm v} 
    = {\bm e}^j \, (\partial_{x_j} \xi_k) \; \;  \partial_{\xi_k} \cdot {\bm v}
    = {\bm r}^k\partial_{\xi_k} \cdot {\bm v},
\end{equation}
where ${\bm e}^j$ is the reciprocal basis of Cartesian coordinates, and ${\bm r}^j$ is the reciprocal basis given by
\begin{equation}
{\bm r}^1 = \frac{{\bm r}_2\times{\bm r}_3}{({\bm r}_2\times{\bm r}_3)\cdot{\bm r}_1} , 
{\bm r}^2 = \frac{{\bm r}_3\times{\bm r}_1}{({\bm r}_3\times{\bm r}_1)\cdot{\bm r}_2} \text{ and }
{\bm r}^3 = \frac{{\bm r}_1\times{\bm r}_2}{({\bm r}_1\times{\bm r}_2)\cdot{\bm r}_3} .
\end{equation}
Note that $({\bm r}_2\times{\bm r}_3)\cdot{\bm r}_1=({\bm r}_3\times{\bm r}_1)\cdot{\bm r}_2=({\bm r}_1\times{\bm r}_2)\cdot{\bm r}_3=V$.
When we use the reciprocal basis ${\bm r}^j$ for $j=1$, 2, and 3, the divergence ${\cal D}_0$ is given by
\begin{equation}
{\cal D}_0({\bm v}) = \sum_{j=1}^3 {\bm r}^j \cdot {\bm w}_j = \sum_{j=1}^3 \hat{\bm r}^j \cdot\frac{{\bm w}_j}{|{\bm r}_j|} 
      ,
     \label{eq:div_euler_norm}
\end{equation}
where $\hat{\bm r}^j={\bm r}^j/|{\bm r}^j|={\bm r}^j|{\bm r}_j|$ are unit vectors for $j=1$, 2, and 3.

Now we consider the Taylor expansion of a function $f({\bm x})$, which is given by 
\begin{equation}
  f({\bm x}+{\bm r}) =  \sum_{k=0}^{n} \frac{1}{k!} ({\bm r}\cdot\nabla)^k f({\bm x}) \, 
  + \, {\cal O}(|{\bm r}|^{n+1})
  .
\end{equation}
Applying the Taylor expansion to the velocity field up to third order, we obtain
\begin{equation}
  \frac{{\bm v}({\bm x}+{\bm r}) - {\bm v}({\bm x})}{|{\bm r}|}
  = (\hat{\bm r}\cdot\nabla){\bm v}({\bm x)}
  + \frac{|{\bm r}|}{2!}(\hat{\bm r}\cdot\nabla)^2{\bm v}({\bm x)}
  + {\cal O}(|{\bm r}|^2)
  .
  \label{eq:Taylor}
\end{equation}
By applying equation~(\ref{eq:Taylor}) to ${\bm w}_j/|{\bm r}_j|$ ($j=1$, 2, and 3) in  equation~(\ref{eq:div_euler_norm}), we get 
%
%
\begin{equation}
{\cal D}_0({\bm v})=
      \sum_{j=1}^3  \hat{\bm r}^j \cdot (\hat{\bm r}_j\cdot\nabla){\bm v}({\bm x}_0)
     + \sum_{j=1}^3  \hat{\bm r}^j \cdot \frac{|{\bm r}_j|}{2!}(\hat{\bm r}_j\cdot\nabla)^2{\bm v}({\bm x}_0)
    + {\cal O}(|{\bm r}|^2)  
  .
\end{equation}
Here, we have $(\hat{\bm r}_j\cdot\nabla){\bm v}=\left(\hat{\bm r}_j\cdot{\bm r}^k \partial_{\xi_k}\right){\bm v}=\frac{1}{|{\bm r}_j|}\left(\delta_j^k \partial_{\xi_k} \right){\bm v}=\frac{1}{|{\bm r}_j|} \partial_{\xi_j}{\bm v} $ for $j=1$, 2, and 3.
Thus, the first term on the right-hand side corresponds to the velocity divergence $\nabla\cdot{\bm v}$ at the position ${\bm x}={\bm x}_0$, and we get
\begin{equation}
{\cal D}_0({\bm v})=
      (\nabla\cdot{\bm v})_{{\bm x}={\bm x}_0} 
     + \sum_{j=1}^3  \hat{\bm r}^j \cdot \frac{|{\bm r}_j|}{2!}(\hat{\bm r}_j\cdot\nabla)^2{\bm v}({\bm x}_0)
    + {\cal O}(|{\bm r}|^2)  
  .
\end{equation}
In the present tessellation-based method,
the volume change is 
calculated for the total volume of the tetrahedrons neighboring
$p_0$ 
\begin{equation}
{\cal D}_p({\bm v}) = \frac{\sum_m {\cal D}_{{\rm tetra},m} V_{{\rm tetra},m}}{\sum_m V_{{\rm tetra},m}},
\end{equation}
with ${\cal D}_{{\rm tetra},m}$ being the divergence computed using equation (\ref{eq:div_euler}) at the $m$-th tetrahedron position neighboring $p_0$. 
%
%
%
%
Thus, the Taylor expansion for ${\cal D}_p$ 
yields
\begin{align}
{\cal D}_p({\bm v}) ~ = ~ &(\nabla\cdot{\bm v})_{{\bm x}={\bm x}_0} \nonumber \\
    & 
     +\sum_m W_m {\bm r}^1_m 
      \cdot|{\bm r}_{m,1}|({\bm r}_{m,1}\cdot\nabla)^2{\bm v}({\bm x}_0) \nonumber \\
    &
     +\sum_m W_m {\bm r}^2_m 
      \cdot|{\bm r}_{m,2}|({\bm r}_{m,2}\cdot\nabla)^2{\bm v}({\bm x}_0) \nonumber \\
    &
     +\sum_m W_m {\bm r}^3_m 
      \cdot|{\bm r}_{m,3}|({\bm r}_{m,3}\cdot\nabla)^2{\bm v}({\bm x}_0) 
     + {\cal O}(|{\bm r}|^2)
    ,
    \label{eq:D_p_sum_V_tetra_W}
\end{align}
%
%
%
%
where $W_m = V_{{\rm tetra}, m}/\sum_{k} V_{{\rm tetra},k}$, and ${\bm r}_m^j$ and ${\bm r}_{m,j}$ for  $j=1$, 2, and 3 are ${\bm r}^j$ and ${\bm r}_{j}$ for each $j$ of the $m$-th tetrahedron, respectively.
Equation~(\ref{eq:D_p_sum_V_tetra_W}) indicates that the proposed method is first-order accurate in space when the effect of $\Delta t$ is negligibly small.
The derivative of ${\bm v}$ is dependent on the tessellation.
We can observe numerically that in the case of a Cartesian grid or a hexagonal tiling, 
the first order terms cancel out which allows to obtain second order convergence.

\section{Numerical validation of the method}
\label{sec:numerical_validation}

In this section we assess the reliability of the proposed method. 
First, we show the advantage of the modified Voronoi tessellation compared with the conventional Voronoi tessellation.
Then we numerically assess the reliability of the method and the convergence rate in space and time 
of the modified Voronoi--based method for computing
a derivative and a sum of derivatives, for a given velocity field. 
The spatial scale is defined as the mean particle distance, $\delta = 2 \pi / N^{1/d}$, where $N$ is the particle number and $d=2, 3$ is the domain dimension.
%
%
Moreover, we study the accuracy as a function of the wavenumber and as a function of the number of particles for a random velocity fields with imposed power-law energy spectra encountered in turbulence.
In the different cases, particles are advected using the explicit Euler scheme.

\subsection{Drawback of Voronoi tessellation}

The reliability of the method of the discrete Voronoi--based divergence for randomly distributed particles has been checked in the Appendix~A of \citet{oujia_matsuda_schneider_2020}. 
We observed a strong correlation between the exact divergence and the Voronoi-based divergence and that the error is most important in strain dominated regions. 
The Pearson correlation is a measure of linear correlation between two data sets. In the case where this coefficient does not tend to $1$ no convergence in space or time is achieved. 
We also observed that when the number of particles increases, the Pearson correlation increases but then saturates at a value of $R[{\cal D}({\bm u}_p), (\nabla \cdot {\bm u})_p] = 0.936$, where $R[\cdot,\cdot]$ is the Pearson correlation coefficient between two sampled quantities.
The coefficient does not tend to $R[{\cal D}({\bm u}_p), (\nabla \cdot {\bm u})_p] = 1$ due to a geometrical error.
The reason is that the deformation of a Voronoi cell is not exactly the same as the deformation of a fluid volume in the continuous setting.
One of the geometrical effects can be explained by behavior similar to the leverage effect induced by the motion of circumcenters, as explained with the following example.
\begin{figure}
\centering
\includegraphics[width=0.45\linewidth]{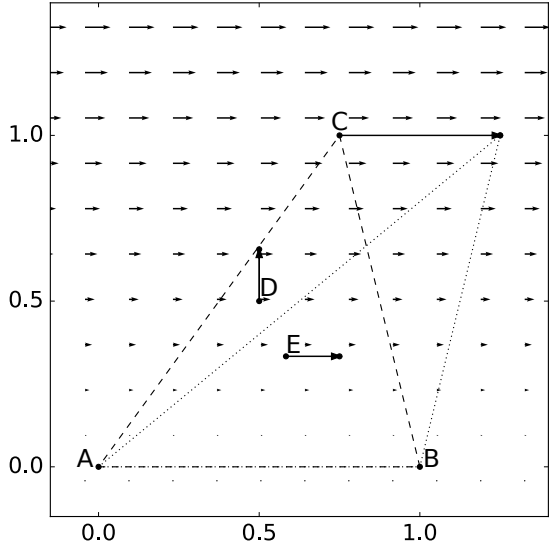}
\caption{Motion of a triangle formed by particles A, B and C, the circumcenter D and the center of gravity E for a linear shear flow. The arrows indicate the particle motion. 
}
\label{fig:Linear_shear_}
\end{figure}
We consider a linear shear flow, ${\bm u}(x,y)=(y,0)$ advecting a triangle of three particles $A$, $B$ and $C$ with $A=(0,0)$, $B=(1,0)$ and $C=(3/4,1)$, represented in figure \ref{fig:Linear_shear_}.
Here $D$ is the center of the circumscribed circle. We observe that when the particles are advected by the velocity ${\bm u}$, the circumcenter $D$ is shifted in $y$ direction, while it should shift towards the same direction as the other points in order to allow the volume of a cell composed by this vertex to be deformed in a similar way as the fluid flow.
Additionally, no angle of the triangle should be near $\pi$ (i.e. so-called flat triangles) because then a large leverage effect is created,
and therefore the variation of the position of the circumcenter is not of the same order of magnitude as that of the triangle.
The motion of the point $D$ can hence result in a variation of the volume of a cell even in the case of zero divergence or in a variation that is not of the right magnitude due to the leverage effect.
%
Moreover, if we reduce the size of the cell and recalling that the divergence is defined as the ratio of the volume variation and 
the volume, the error remains unchanged and persists regardless of the scale.
In contrast, if we define $E$ as the center of gravity of the triangle, we can see that its motion is similar to that of the shear. 
Using the modified Voronoi tessellation for the same velocity field as in Appendix~A of \citet{oujia_matsuda_schneider_2020}, we obtain that the Pearson correlation tends to $R[{\cal D}({\bm u}_p), (\nabla \cdot {\bm u})_p] = 1$.
For this reason, we will focus in the following only on the modified Voronoi tessellation.

\subsection{Computation of a spatial derivative }
\label{sec:Computation_of_a_spatial_derivative}

Subsequently, we verify numerically the validity of the method for the computation of a derivative using the modified Voronoi--based method in two and three dimensions. 
To this end, we consider randomly distributed particles in a $2\pi$-periodic domain $\Omega \subset \mathbb{R}^d$ ($d=2,3$) advected by a given velocity field. The derivative in a given direction is computed using the variation of the modified Voronoi volume. In two dimensions, we compare the modified Voronoi--based method with conventional Eulerian methods.

\begin{figure}
\centering
(a)
\includegraphics[width=0.40\linewidth]{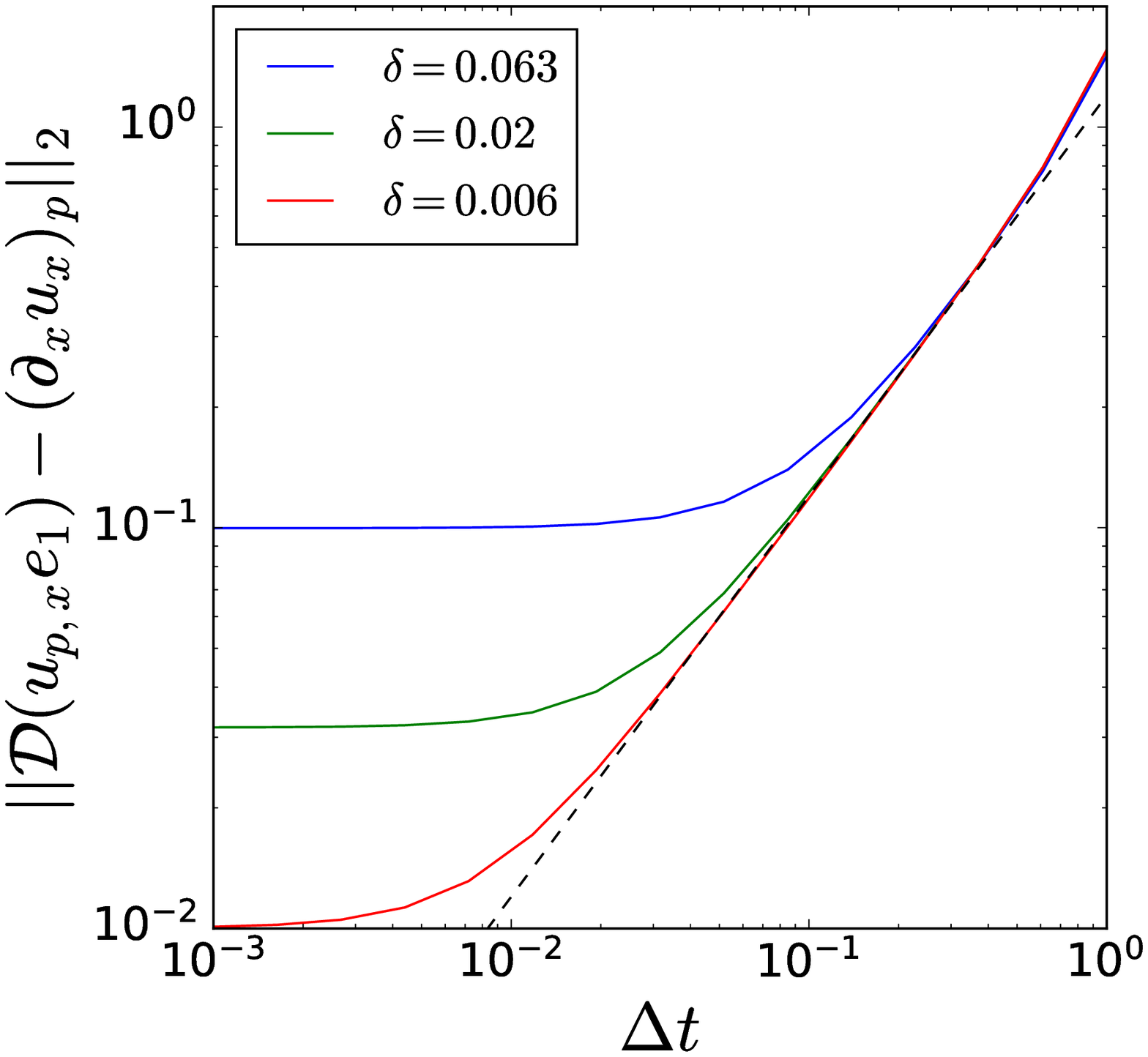}
(b)
\includegraphics[width=0.40\linewidth]{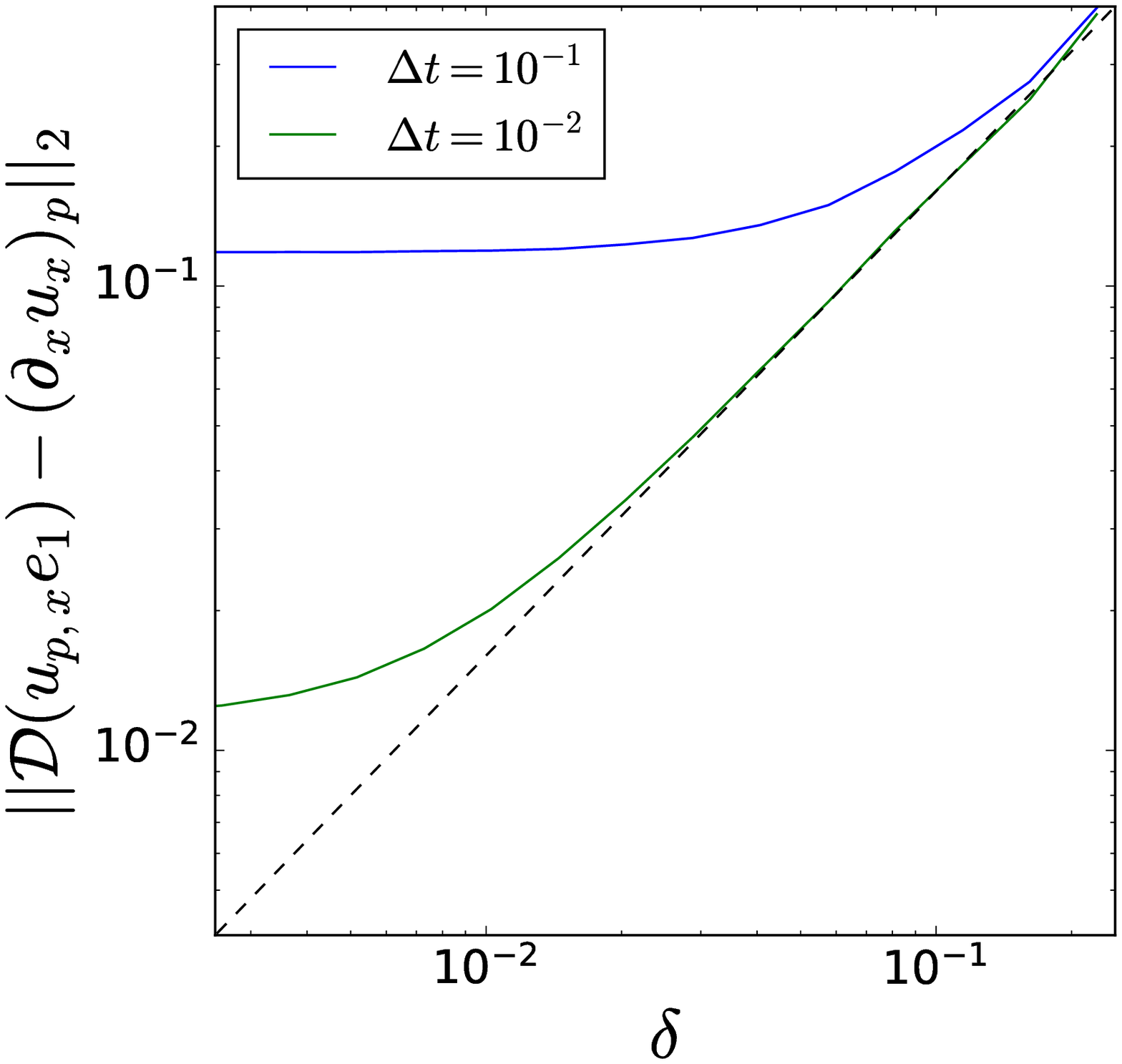}\\
(c)
\includegraphics[width=0.40\linewidth]{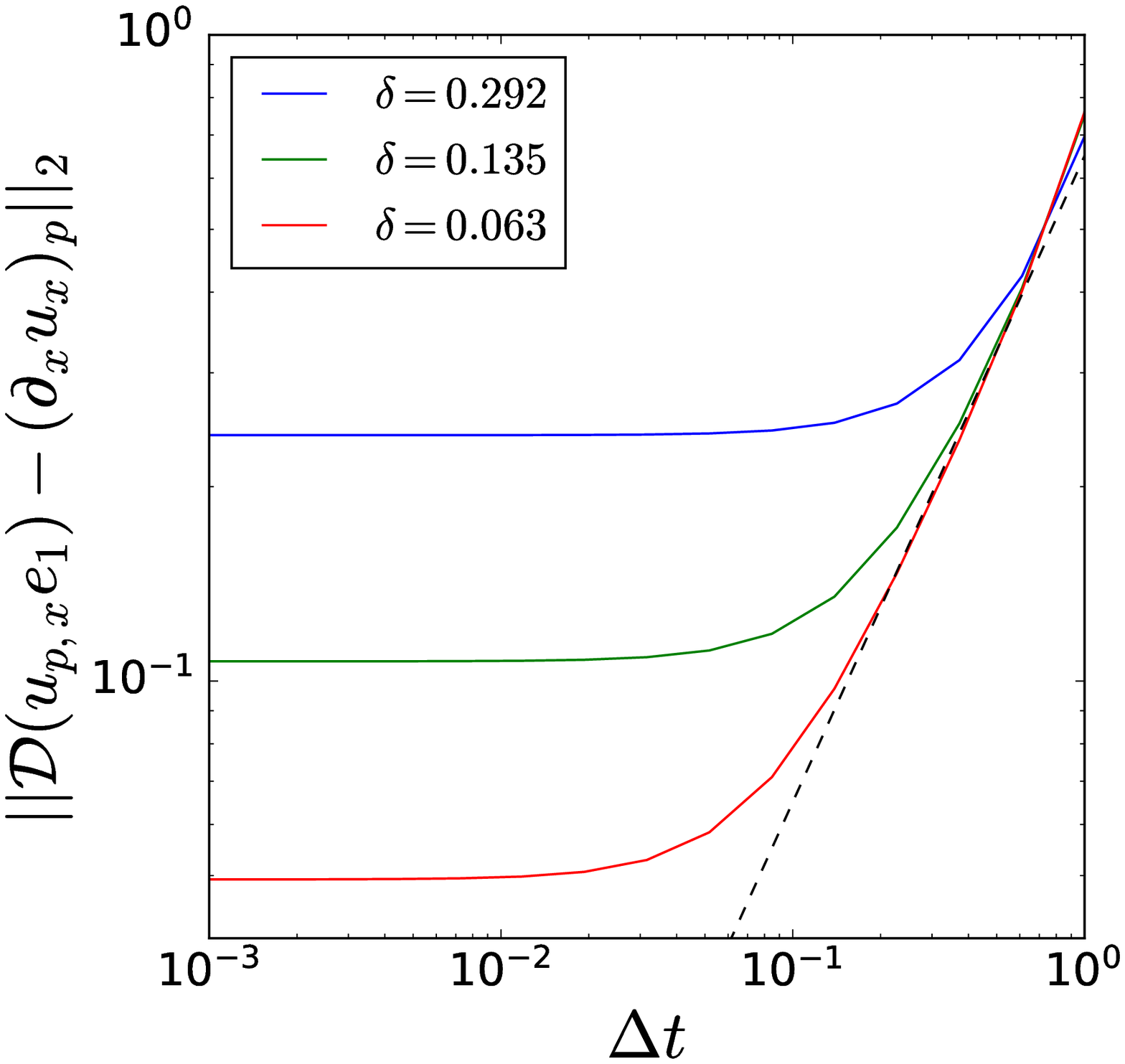}
(d)
\includegraphics[width=0.40\linewidth]{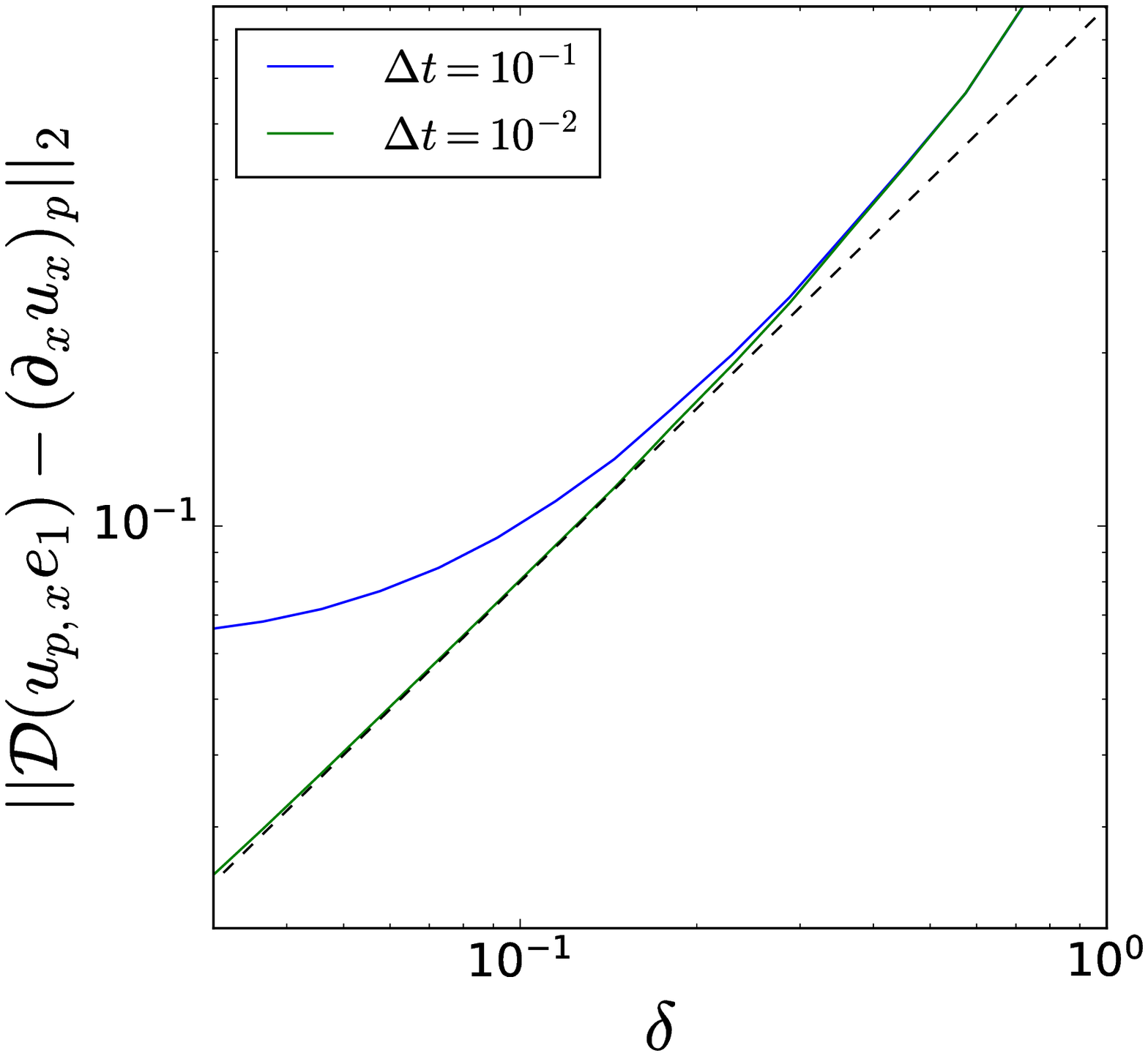}
\caption{ $L^2$ error of the $x$-derivative of the velocity field $\bm u$, given by eqs. (21) and (22) for two and three dimensions respectively, 
as function of $\Delta t$ (a,c) and $\delta$ (b,d) for two (a,b) and three dimensions (c,d) for randomly distributed particles using the modified Voronoi--based method. Dashed lines indicate slope 1. 
}
\label{fig:Error_Derivative_2D}
\label{fig:Error_Derivative_3D}
\end{figure}

\subsubsection{Two dimensions}
\label{sec:Computation_of_a_spatial_derivative_2D}

In order to check the reliability of the method in two dimensions in terms of space and time convergence, we consider the following velocity field 
%
${\bm u} = (u_x, 0)^T$, where $u_x = \cos(x)\cos(y)$ 
whose partial derivative in $x$-direction is, $\partial_x u_x  = -\sin(x)\cos(y)$. 
%
Figure~\ref{fig:Error_Derivative_2D}(a) and (b) show the $L^2$ error 
of $\displaystyle \partial_x { u}_x $ for two-dimensional randomly distributed particles as a function of $\Delta t$ for a fixed value of the mean particle distance $\delta$, and as a function of $\delta$ for fixed $\Delta t$, respectively. We observe first order convergence in both time and space. 
The plateau is caused by the interdependence between $\Delta t$ and $\delta$, which is implied by the Lagrangian aspect of the method.

\begin{figure}
\centering
(a)
\includegraphics[width=0.45\linewidth]{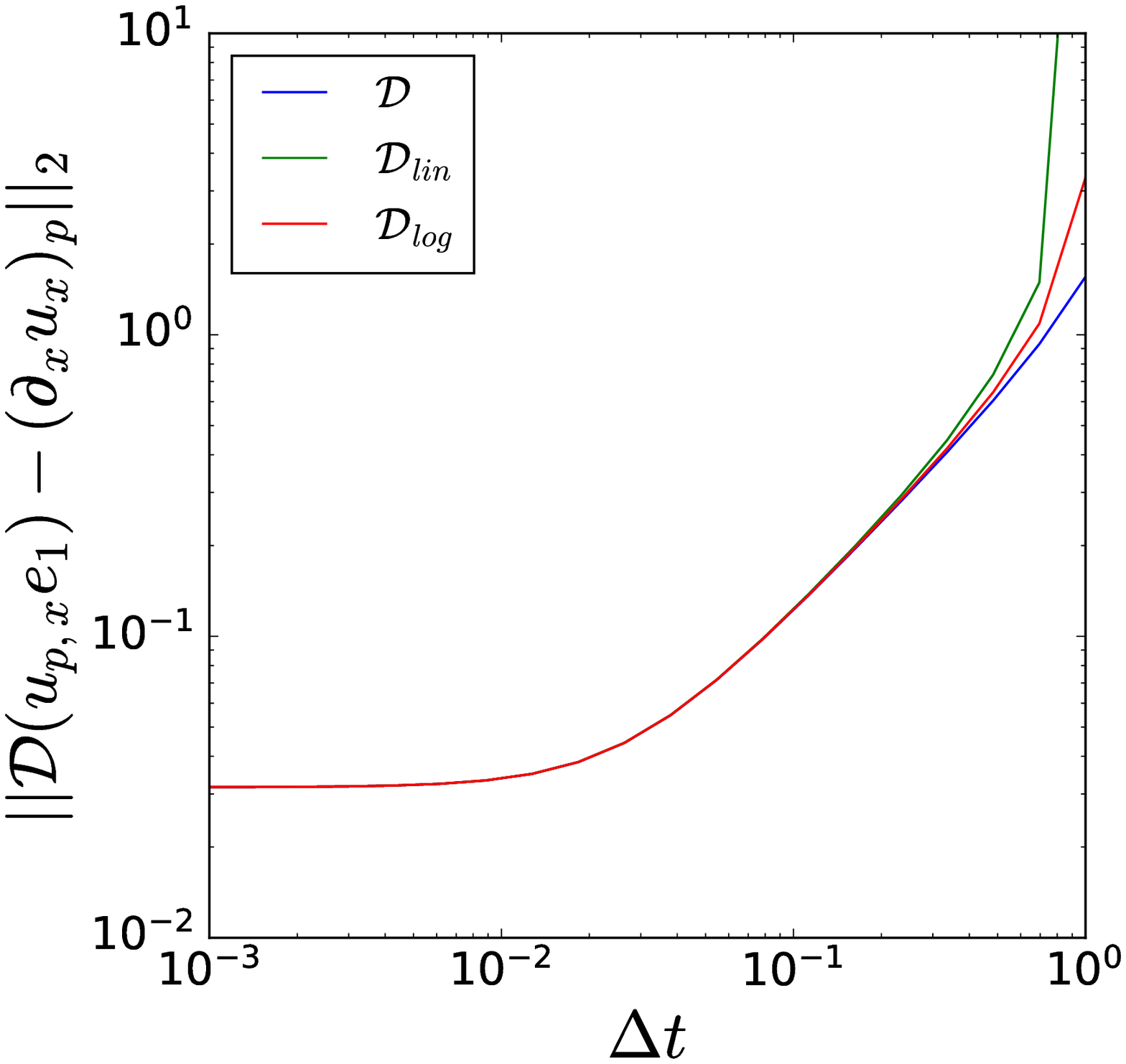}
(b)
\includegraphics[width=0.45\linewidth]{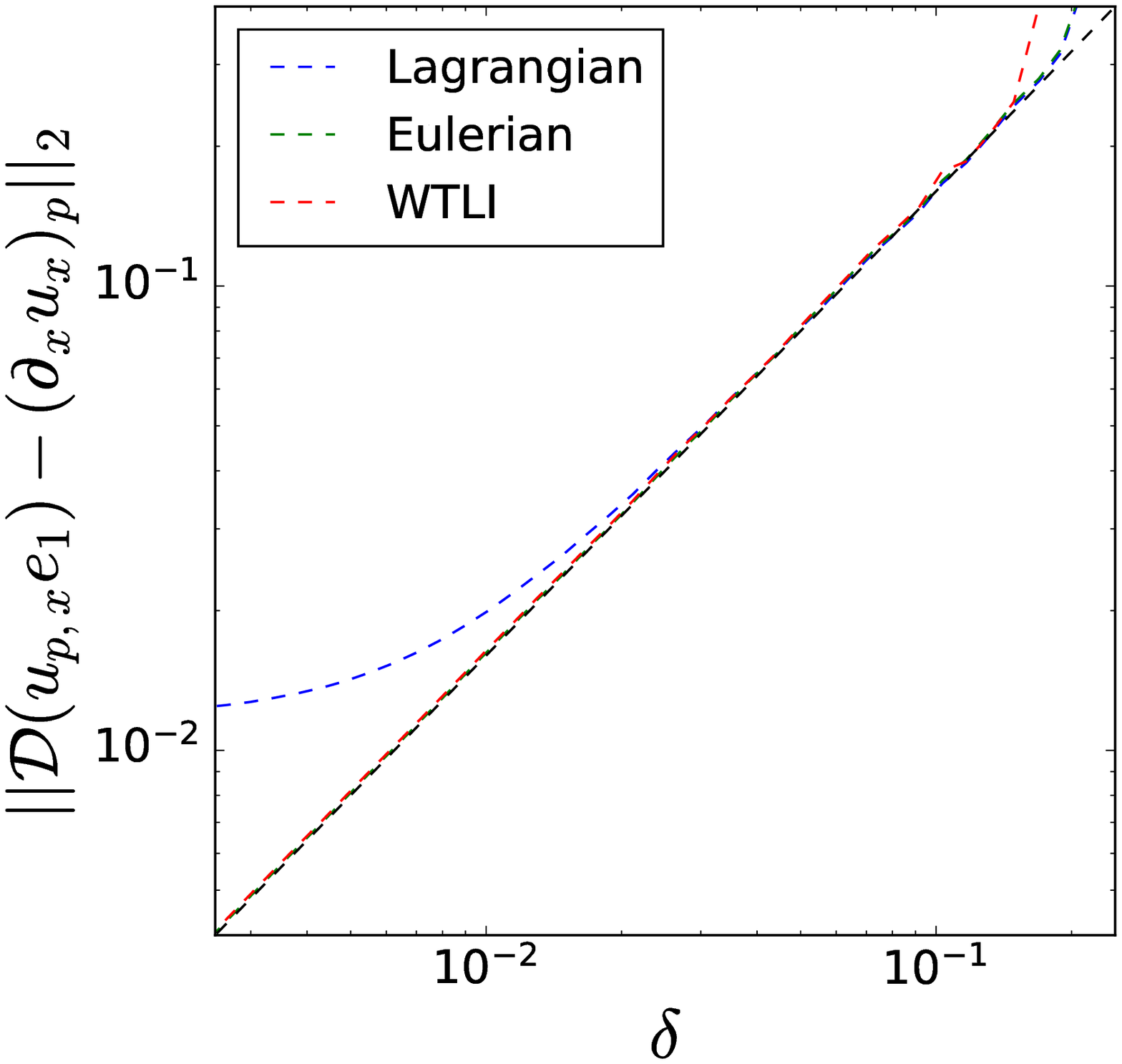}
\caption{
$L^2$ error of the gradient as function of $\Delta t$ (a) in two dimensions for randomly distributed particles using the modified Voronoi--based method with different temporal approximation of the particle density, 
i.e., $\cal D$, ${\cal D}_{lin}$, and ${\cal D}_{\log}$ defined as eqs.~(\ref{eq:div_density_tessellation}), (\ref{eq:div_density_tessellation_lin}), and (~\ref{eq:div_density_tessellation_log}), respectively.
$L^2$ error of the $x$-derivative of the velocity field $\bm u$, given by eq. (21), 
using Lagrangian and Eulerian modified Voronoi--based method and WTLI method (b) computed for two-dimensional randomly distributed particles as function of $\delta$. The dashed line has a slope $1$.
} 
\label{fig:Error_Dynamique_Static}
\end{figure}

Figure \ref{fig:Error_Dynamique_Static}(a) shows the $L^2$ error of $\displaystyle \partial_x {u}_x$ for two-dimensional randomly distributed particles as a function of $\Delta t$ for a fixed number of particles $N=10^5$ using different temporal approximation of the particle density defined in subsection \ref{sec:div} by equation~(\ref{eq:div_density_tessellation}), equation~(\ref{eq:div_density_tessellation_lin}) and equation~(\ref{eq:div_density_tessellation_log}).
We can observe that, as expected, at low $\Delta t$ the error is similar, and we discern first order decay and then saturation at $3 \times 10^{-2}$ due to the spatial error. 
However, when $\Delta t$ is large, the error is the lowest for equation~(\ref{eq:div_density_tessellation}), highest for equation~(\ref{eq:div_density_tessellation_lin}) and in between for equation~(\ref{eq:div_density_tessellation_log}). 
Thus, the divergence based on equation~(\ref{eq:div_density_tessellation}) is the best approximation among those based on the three proposed equations.

Now we compare these results with an existing gradient computation
method described in \citet{sozer2014gradient} which utilizes the Green-Gauss Gradient Method. 
The Weighted Tri-Linear Face Interpolation (WTLI) method from is used for simplicity of implementation, 
since the other methods (e.g., Weighted Least Squares Face Interpolation or Least Squares Gradient Method) presented in \citet{sozer2014gradient} give
similar results.
The gradient of ${\bm u}$ in a closed area $V$ is defined as
\begin{equation} 
\nabla {\bm u}_p = \frac{1}{V_p}\sum^{N_{\rm edges}}_{e=1} \overline{{\bm u}}_e \hat{\bm n}_e A_e 
\end{equation}
where $N_{\rm edges}$ is the number of edges of the cell formed by the tessellation, $A_e$ is the length of the edge $e$,
$\hat{\bm n}_e$ is the edge unit normal vector and $\overline{{\bm u}}_e$ is the average of ${\bm u}$ over the edge $e$.
Using the WTLI method, the value of $\overline{{\bm u}}_e$ is given by a barycentric linear interpolation of the velocity at the midpoint of the edge $e$ using the velocity of three particles 
adjacent to the center of the edge $e$.

Figure \ref{fig:Error_Dynamique_Static} (b) shows the $L^2$ error of $\displaystyle \partial_x { u}_x$ using the Lagrangian ($\Delta t = 10^{-2}$) and Eulerian modified Voronoi--based method, i.e., computed respectively with equation (\ref{eq:div_density_tessellation}) and (\ref{eq:div_euler}), and the WTLI method computed for two-dimensional randomly distributed particles as function of $\delta$ using the exact value of the particle velocity ${\bm u}$. 
We observe that the results for the different methods are similar.
For the different methods, we find first order in space, and a saturation of the error for the Lagrangian method due to the 
interdependence with $\Delta t$.
The results obtained with the modified Voronoi based method and with the WTLI method are thus comparable. 
Hence the Eulerian modified Voronoi-based method allows to overcome the influence of $\Delta t$, compared to the Lagrangian version. However, it is observed that our implementation of the Lagrangian method is faster than the one for the Eulerian modified Voronoi-based method.
For this reason the Lagrangian method is a new efficient alternative to the existing Eulerian one.
%
%

\subsubsection{Three dimensions} 

Now we examine the space and time convergence in three dimensions and we choose the following velocity field 
${\bm u} = (u_x, 0, 0)^T$, where $u_x = \sin(x)\cos(y)\cos(z)$ 
and whose partial derivative in $x$-direction is $\partial_x u_x = \cos(x)\cos(y)\cos(z)$.
Figure~\ref{fig:Error_Derivative_3D}(c) shows the $L^2$ error of $\displaystyle \partial_x { u}_x$ for three-dimensional randomly distributed particles as a function of $\Delta t$ for a fixed value of $\delta$, 
and figure~\ref{fig:Error_Derivative_3D}(d) as a function of $\delta$ for a fixed $\Delta t$. We can observe first order in time  
and space similar to what has been observed for the two dimensional case.
Again the saturation of the error for large values of $\delta$ and small $\Delta t$ is due to the 
interdependence 
between $\Delta t$ and $\delta$.
%
Therefore, we can conclude that the method converges with first order in time and in each space direction for randomly distributed particles, as expected and shown theoretically in section~\ref{sec:euler}. 


\subsection{Computation of the sum of derivatives}

Now we numerically verify the validity of the method for the computation of a sum of derivatives required for the divergence or the curl of the particle velocity. Indeed, the method being initially designed for the computation of the divergence, gives us the ability to compute simultaneously the sum of derivatives of a function by rearranging the different components of the function. In this section we will compute the divergence and show that the method converges, as expected, with the same order in space and time for the sum of partial derivatives.

\begin{figure}
\centering
(a)
\includegraphics[width=0.40\linewidth]{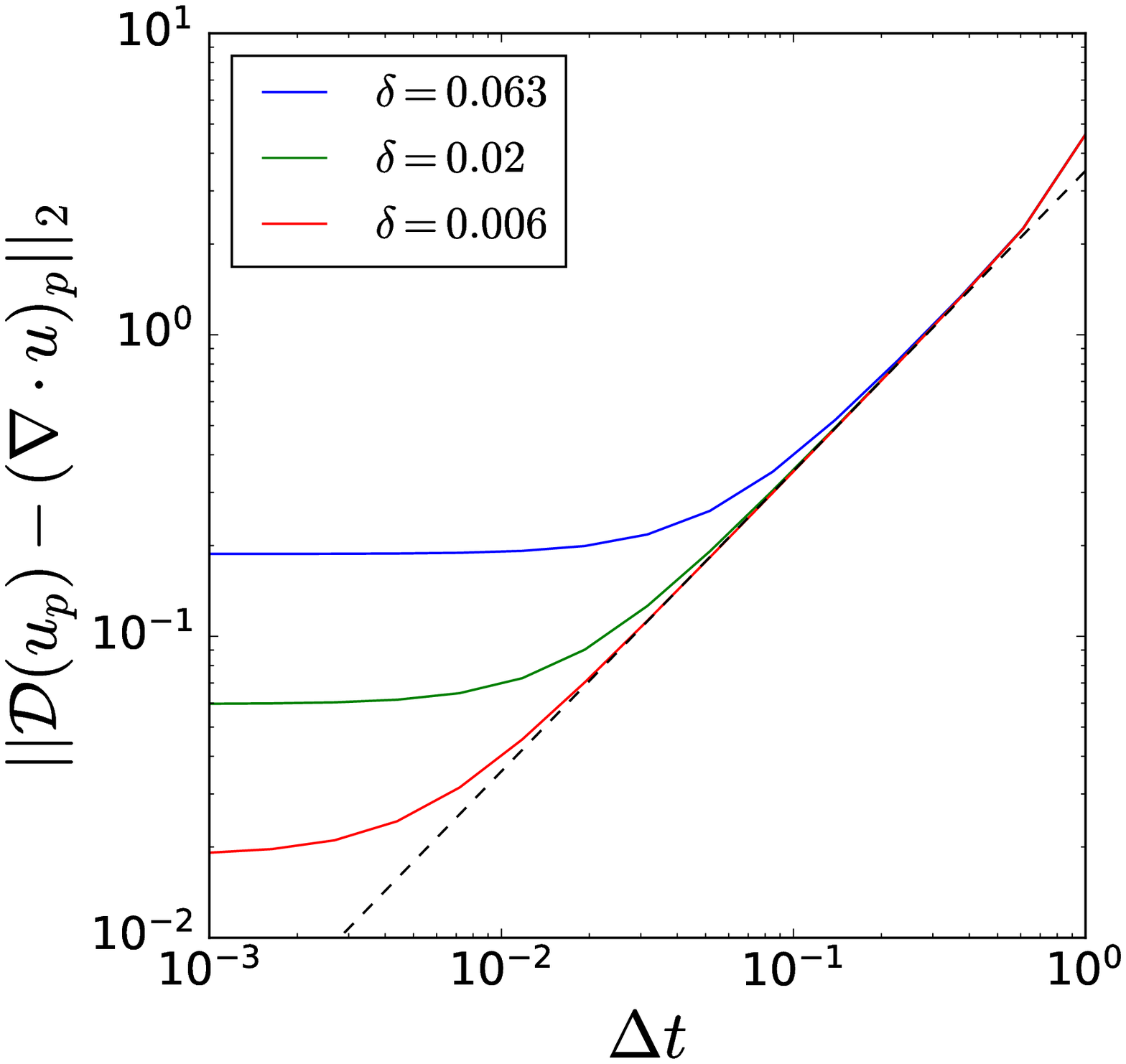}
(b)
\includegraphics[width=0.40\linewidth]{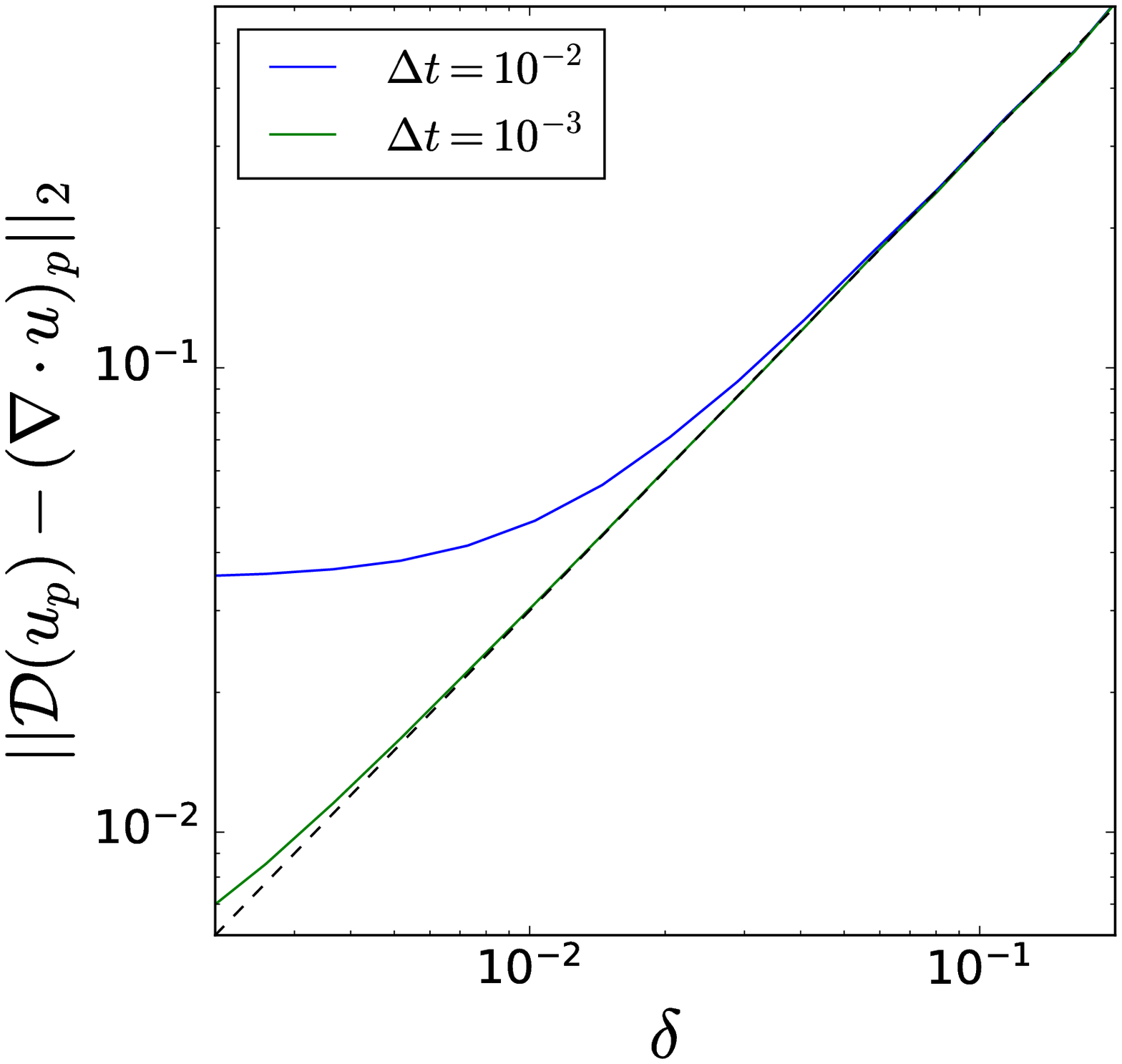}\\
(c)
\includegraphics[width=0.40\linewidth]{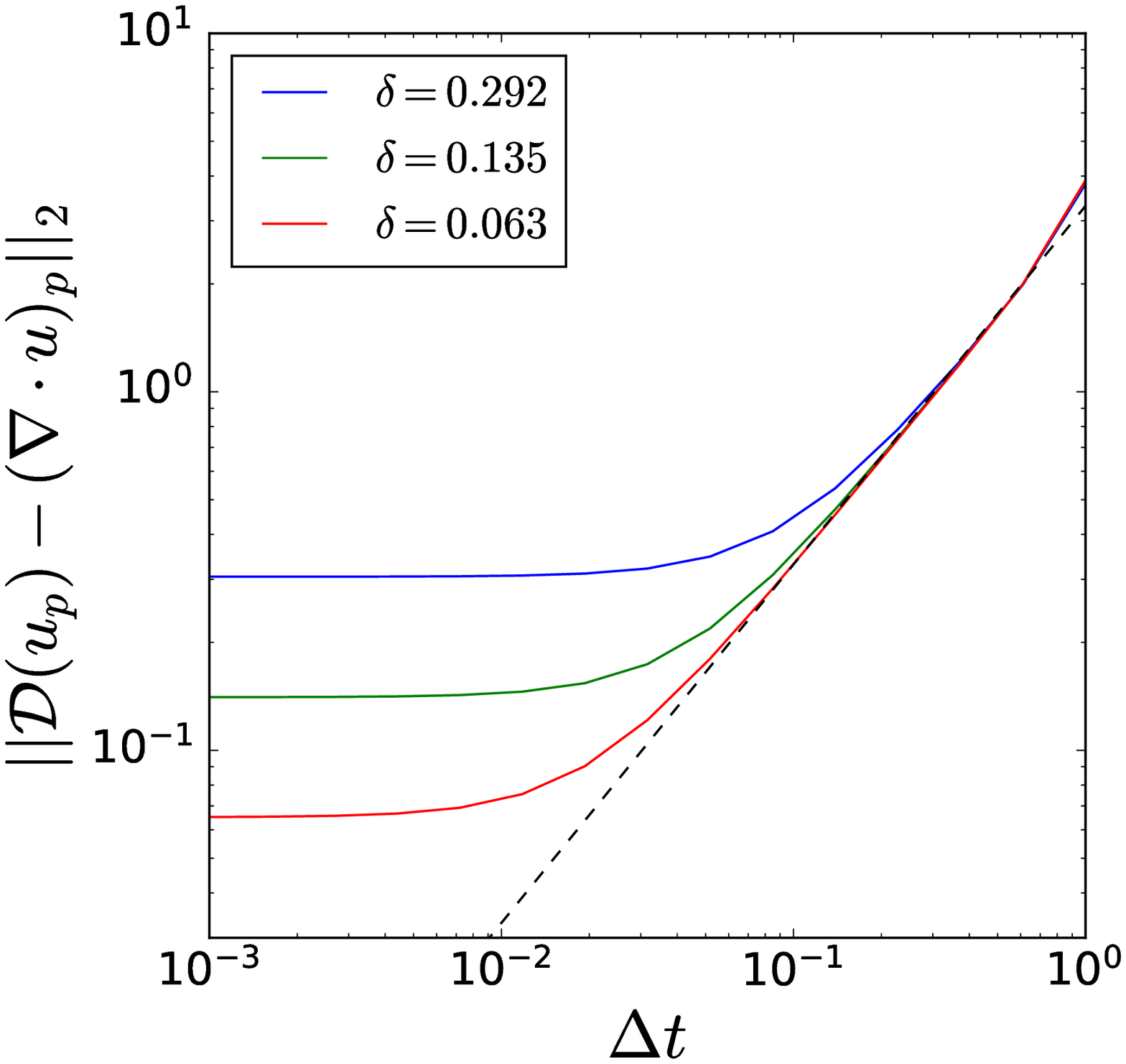}
(d)
\includegraphics[width=0.40\linewidth]{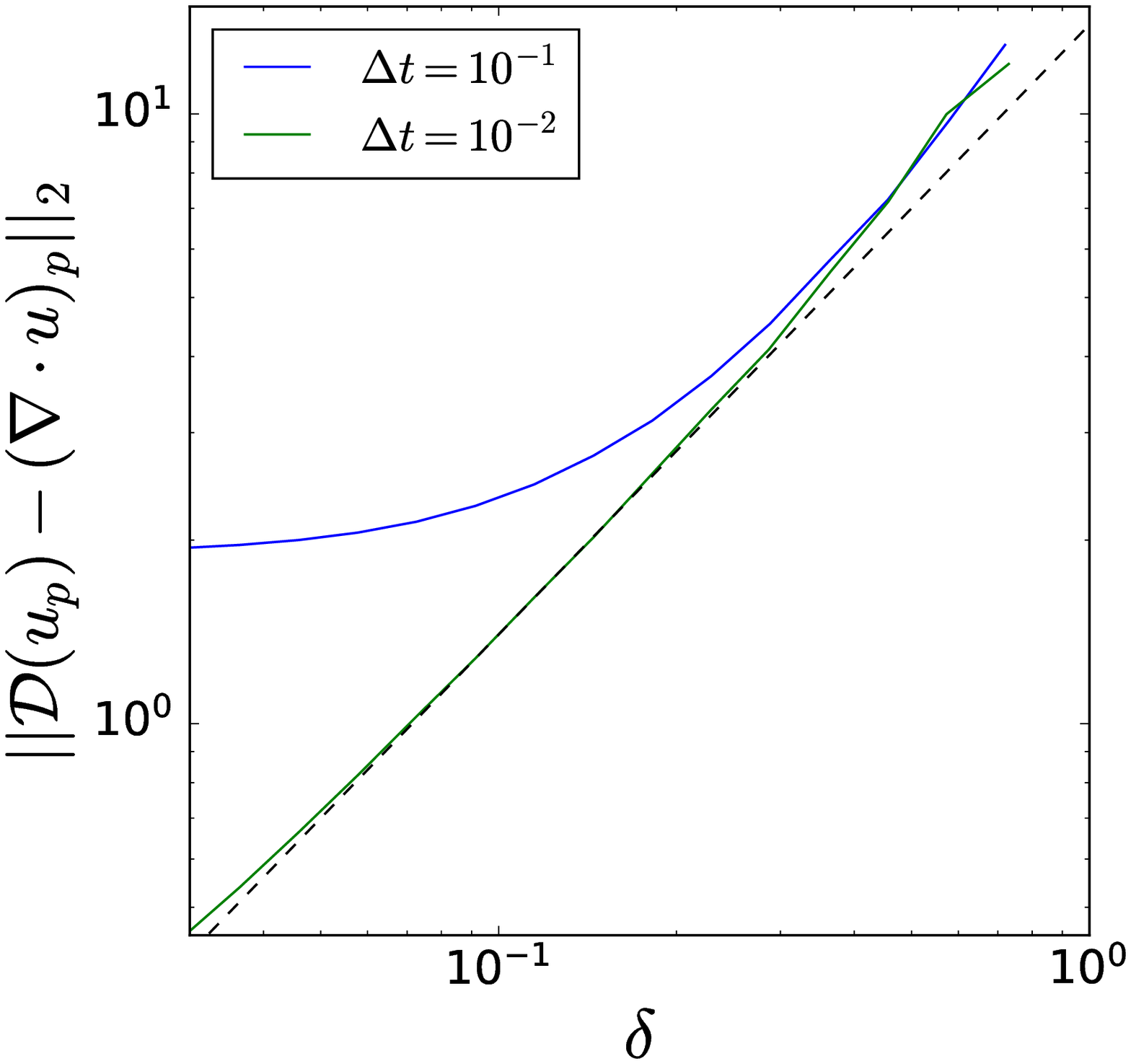}
\caption{ $L^2$ error of the divergence as a function of $\Delta t$ (a,c) and $\delta$ (b,d) for two (a,b) and three dimensions (c,d) for randomly distributed particles using the modified Voronoi--based method. Dashed lines indicate slope 1. 
}
\label{fig:erreur_del_curl_2D}
\label{fig:erreur_del_curl_3D}
\end{figure}

We analyze the convergence of the discrete divergence operator in space and time considering the velocity field 
\begin{equation}
    ~~~~~~~~~~~~~~~~~~~~~~~~~~~~~~~~~~~~~~~~~~~~~~~
    {\bm u}(x, y) = \begin{pmatrix}
    ~~\cos(x)\cos(y)\\ -\sin(x)\sin(y)
    \end{pmatrix}
    ~~~~~~~~~~~~~~~~~~~~~~~~~~~~ (x,y) \in \Omega
\end{equation}
in two dimensions, whose divergence is  $\nabla \cdot {\bm u}(x, y) = -2\sin(x)\cos(y)$, and in three dimensions
\begin{equation}
    ~~~~~~~~~~~~~~~~~~~~~~~~~~~~~~~~
    {\bm u}(x, y, z) = \begin{pmatrix}
    \sin(x)\cos(y)\cos(z)\\ \cos(x)\sin(y)\cos(z)\\ \cos(x)\cos(y)\sin(z),
    \end{pmatrix}
    ~~~~~~~~~~~~~~~~ 
    (x,y,z) \in \Omega
\end{equation}
whose divergence is $\nabla \cdot {\bm u}(x, y, z) = 3\cos(x)\cos(y)\cos(z)$.
Figure \ref{fig:erreur_del_curl_2D} shows the $L^2$ error $\Vert {\cal D}({\bm u}_p) - (\nabla \cdot {\bm u})_p \Vert_2$ as a function of $\Delta t$ for a fixed value $\delta$ (a, c) and as a function of $\delta$ for fixed $\Delta t$ (b, d) in two (a,b) and three dimensions (c,d). 
As expected, the results are comparable to those in figure \ref{fig:Error_Derivative_2D}, i.e. we observe a decay of the $L^2$ error with slope $1$ for both, the error in time and space  in two and three dimensions. This confirms that the method allows to compute simultaneously the sum of derivatives and thus the divergence or curl with the same order of approximation.

A similar analysis for the curl of a vector field composed of sine and cosine functions has been performed (not shown here), and alike results with the same order of convergence are found.

\subsection{Curl of synthetic turbulent velocity fields and required particle number}

We consider synthetic turbulent velocity fields, i.e., random fields with a given correlation, for advecting randomly distributed particles. The energy spectrum satisfies a power law, with slopes typically observed in two- and three-dimensional turbulence.
The aim is to determine the required number of particles as a function of the wavenumber, in order to achieve a strong correlation between the exact values and those computed using the tessellation-based method.
%
%
The derivative of a sine wave 
with given wavenumber $k$ in a $2\pi$--periodic domain is computed
for ${\bm u}^{2D}_k(x, y) = (\sin(kx),0)$ in two dimensions and ${\bm u}^{3D}_k(x, y, z) = (\sin(kx),0,0)$ in three dimensions 
for different $k \delta$ values.
\begin{figure}
\centering
(a)
\includegraphics[width=0.45\linewidth]{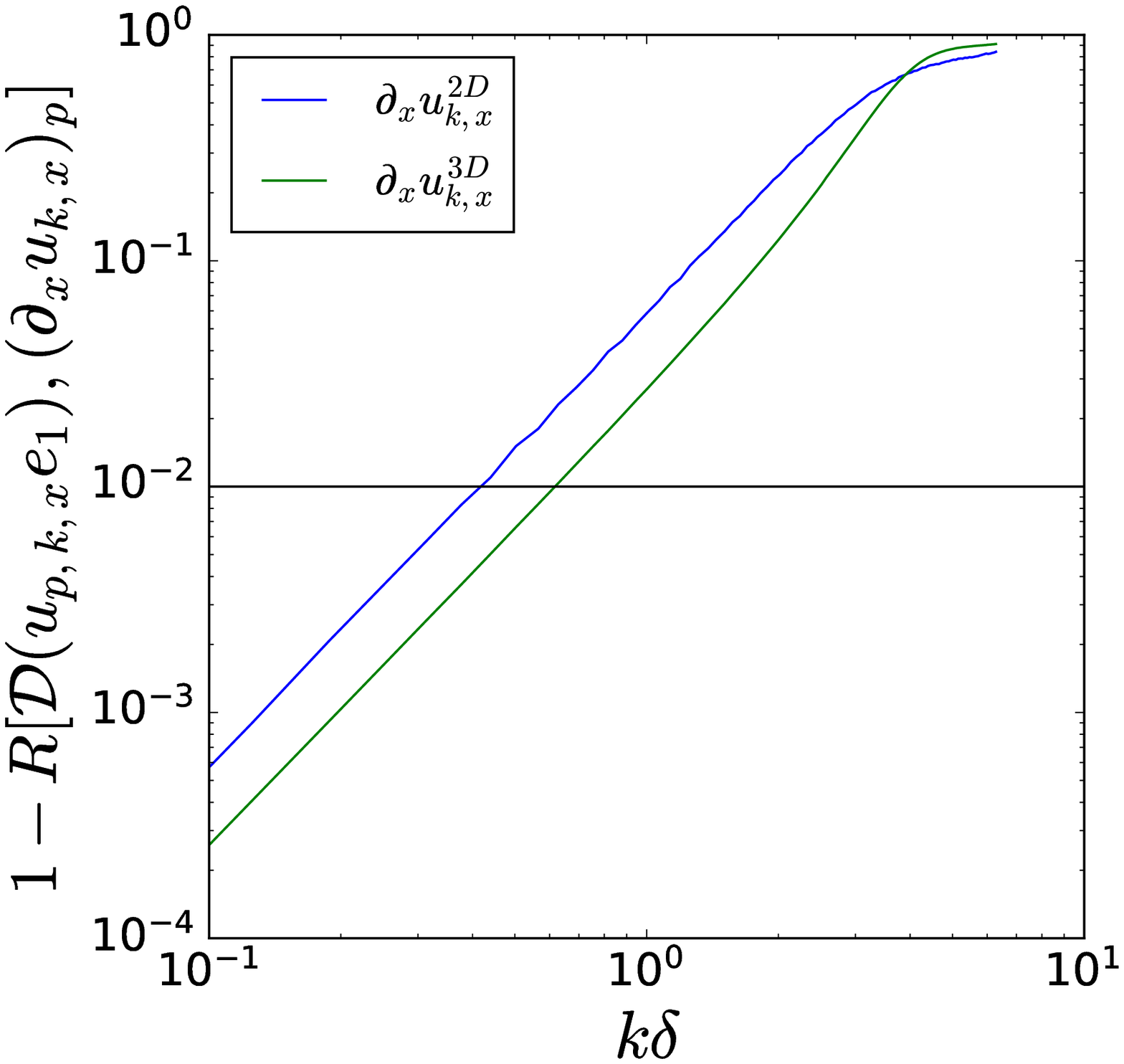}
(b)
\includegraphics[width=0.45\linewidth]{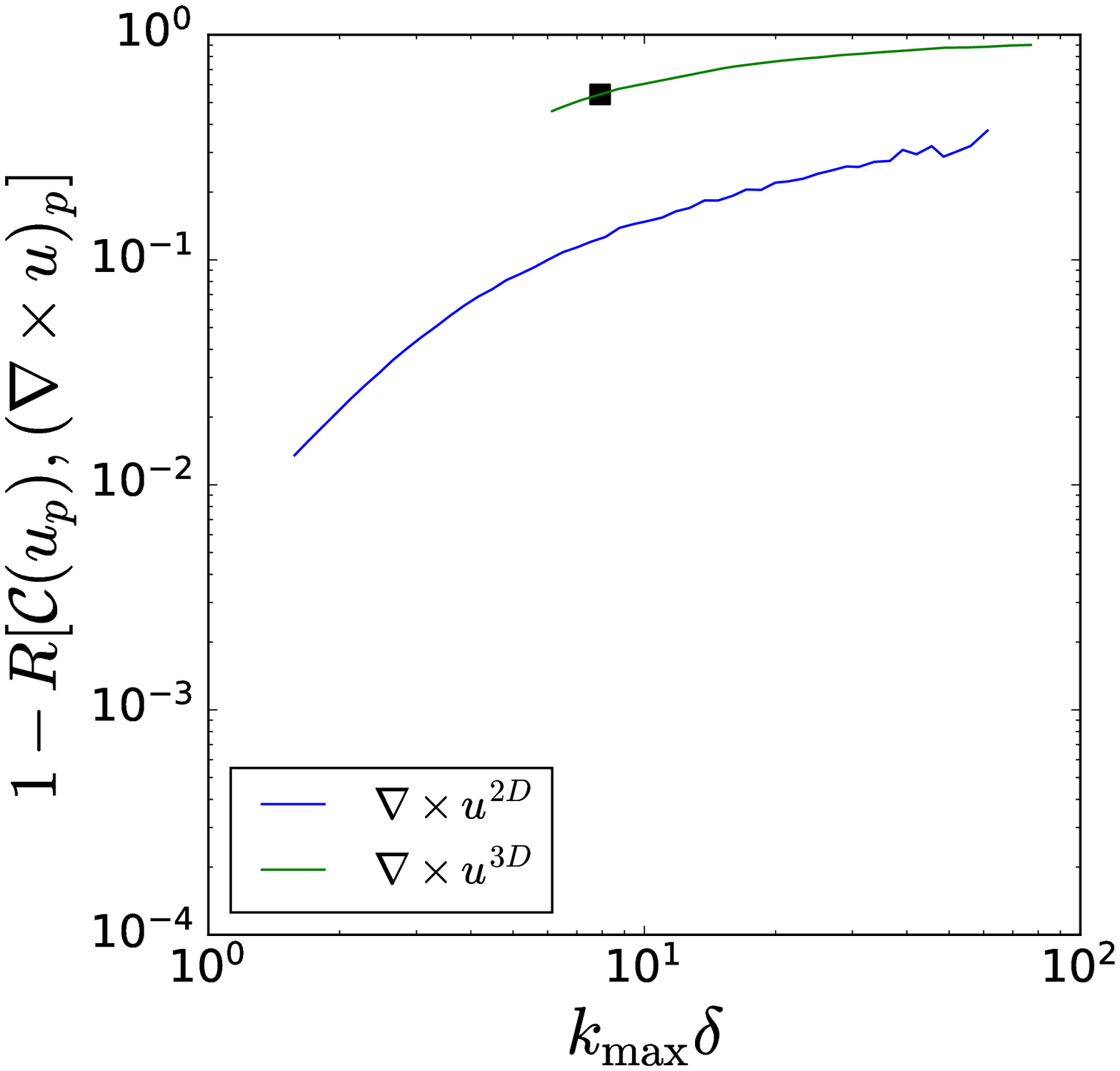} 
\caption{
Pearson correlation coefficient $R$ 
between the exact derivative values and computed values using the modified Voronoi--based derivative as function $k \delta$ (a) in two and three dimensions.  
Pearson correlation coefficient $R$ between the exact curl values and computed values using the modified Voronoi--based curl as function the number of particles for velocity fields with imposed power-law energy spectra (b) in two and three dimensions.
The black horizontal line indicates the Pearson correlation value $0.99$.
The black square indicates the theoretical lower bound for the analyzed DNS data, for details we refer to section~\ref{sec:application_to_fluid_particles}.
} 
\label{fig:pearson_2D_3D_frequence}
\end{figure}
Then the Pearson correlation between the exact derivative value and the numerical one is calculated and shown in Figure~\ref{fig:pearson_2D_3D_frequence} as function of $k \delta$.
As expected, the smaller the value of $k \delta$, the higher the Pearson correlation. 
%
Thus, to have a correlation of $0.99$ an approximate number of particles of $N = (2\pi k)^2/(4.177\times10^{-1})^2$ is needed in two dimensions, and $N = (2\pi k)^3/(6.173\times10^{-1})^3$ in three dimensions. 
The values of $4.18\times10^{-1}$ and $6.17\times10^{-1}$ are given by the intersection of the black horizontal line with the different curves in figure~\ref{fig:pearson_2D_3D_frequence}.
%
%

%
In the next step, 
we apply the method to synthetic turbulence whose energy spectrum satisfies a power law behavior.
%
%
In order to have a lower bound estimation of the accuracy of the method in the case of the computation of the vorticity, i.e. the curl of the velocity of a turbulent flow, we consider a realization of random functions ${\bm u}^{2D}$ and ${\bm u}^{3D}$, which are the sum of sines with a random phase defined as, 
\begin{equation}
{\bm u}^{2D}(x, y) = 
    \sum_{k=1}^{k_{\rm max}}
    {E(k)^{1/2}}
    \begin{pmatrix}
        \sin(kx+r_{k,0}) + \sin(ky+r_{k,1}) \\
        \sin(kx+r_{k,2}) + \sin(ky+r_{k,3}) 
    \end{pmatrix}
\quad (x,y) \in \Omega
\end{equation}
\begin{equation}
{\bm u}^{3D}(x, y, z) = 
    \sum_{k=1}^{k_{\rm max}}
    {E(k)^{1/2}}
    \begin{pmatrix}
        \sin(kx+r_{k,0}) + \sin(ky+r_{k,1}) + \sin(kz+r_{k,2})\\
        \sin(kx+r_{k,3}) + \sin(ky+r_{k,4}) + \sin(kz+r_{k,5})\\
        \sin(kx+r_{k,6}) + \sin(ky+r_{k,7}) + \sin(kz+r_{k,8})
    \end{pmatrix}
    \quad (x,y, z) \in \Omega
\end{equation}
where $r_{k,j}\in[0,2\pi[$ for $j=0,\cdots,8$ 
are uniformly distributed random numbers, and $E(k) = k^{-3}$ in two dimensions and $k^{-5/3}$ in three dimensions, which are the theoretical prediction of the spectrum in each of the dimensions for a turbulent flow.
We take $k_{\rm max}=256$ to be able to compare the results with a numerical simulation with a resolution of $512^{d}$ ($d=2,3$).
%
Figure \ref{fig:pearson_2D_3D_frequence}(b) shows the Pearson correlation $R[\mathbfcal{C}({\bm u}_p^{2D}),(\nabla \times {\bm u}^{2D})_p]$ and $R[\mathbfcal{C}({\bm u}_p^{3D}),(\nabla \times {\bm u}^{3D})_p]$ as function of $k_{\rm max} \delta$. 
%
For two and three dimensions, we observe that the correlation values in figure \ref{fig:pearson_2D_3D_frequence}(b) are higher than those in figure \ref{fig:pearson_2D_3D_frequence}(a). This difference is especially noticeable at $k_{max}\delta = 6$, which is the only common data point for both curves.
%
This can be explained by the fact that the energy spectrum $E(k)$ rapidly decays with $k$ and thus the high wavenumber coefficients in the velocity fields ${\bm u}^{2D}$ and ${\bm u}^{3D}$ have low amplitude. 
\section{Application to turbulence}
\label{sec:application_to_turbulence}

In this section we test our method for data from state-of-the art Direct Numerical Simulation (DNS) of particle-laden homogeneous isotropic turbulence (HIT), detailed in \citet{matsuda2021scale}. 
The incompressible Navier--Stokes equations are solved in a $2 \pi$-periodic cube with a fourth-order finite-difference scheme. A large scale forcing is applied to obtain a statistically stationary turbulent flow.
High resolution DNS computations with $N_g^3 = 512^3$ grid points are performed for the Taylor-microscale Reynolds number $Re_\lambda = 204$.
%
Uniformly distributed discrete particles are injected into the fully developed turbulent flow, where the number of fluid and inertial particles is $N_p = 1.07 \times 10^9$.

\begin{figure}
\centering
\includegraphics[width=0.95\linewidth]{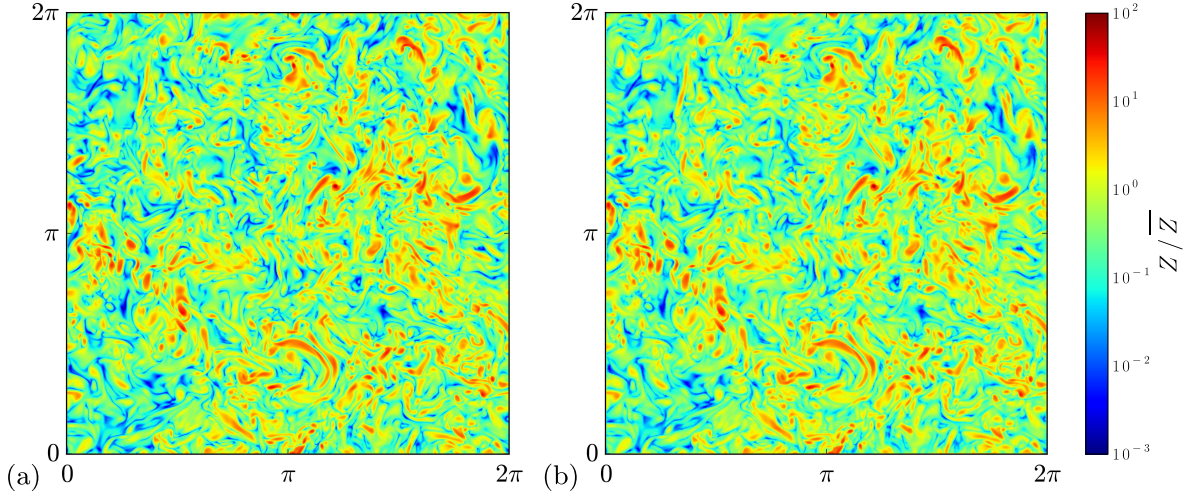}
\caption{
Enstrophy field $Z$ normalized by the mean value of the turbulent flow in a slice of $2\pi/512$ of the domain computed on an Eulerian grid (a) and with the modified Voronoi--based curl (b) where the enstrophy is the $L^2$ norm squared of the curl.
}
\label{fig:Vor_Grid_enstrophy}
\end{figure}


\subsection{Application to fluid particles}
\label{sec:application_to_fluid_particles}
First we consider fluid particles, 
i.e. particles without inertia which follow the streaklines of the flow. 
To get the velocity of fluid particles we sample the fluid velocity at random positions.
%
%
%

%
The curl of fluid particles is computed using the total number for particles advected by the fluid velocity at particle positions.
%
Figure~\ref{fig:Vor_Grid_enstrophy}(a) shows 
the enstrophy $Z=||\nabla\times{\bm u}||_2^2$
in a two-dimensional slice at $z\in[0,2\pi/512]$ 
computed on an Eulerian grid, 
and figure~\ref{fig:Vor_Grid_enstrophy}(b) the enstrophy computed from the curl of fluid velocity by using the modified Voronoi--based method, i.e., $Z_p = ||{\mathbfcal C}({\bm u}_p)||_2^2$.
The square of the curl is normalized by the mean value.
Note that, in figure \ref{fig:Vor_Grid_enstrophy}(b), $Z_p$ is spatially averaged in each box on the equidistant grids with the size of $512^3$. 
%
A perfect agreement is observed between the two methods.
This visual result gives some confidence and good expectations applying the modified Voronoi--based curl to real flow data either obtained by DNS or by particle tracking velocimetry (PTV) of experiments.


\begin{figure}
\centering
(a)
\includegraphics[width=0.45\linewidth]{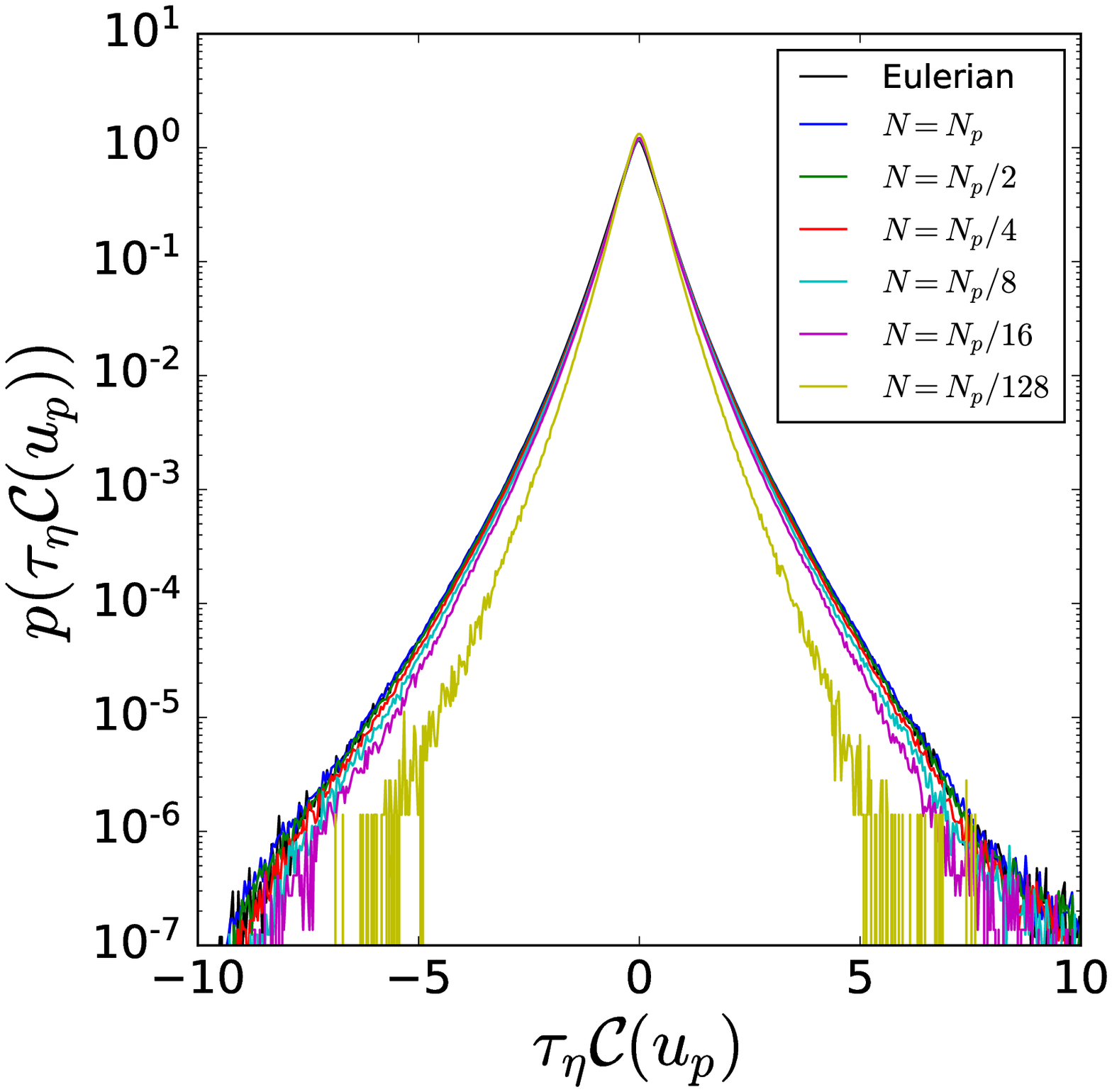}
(b)
\includegraphics[width=0.45\linewidth]{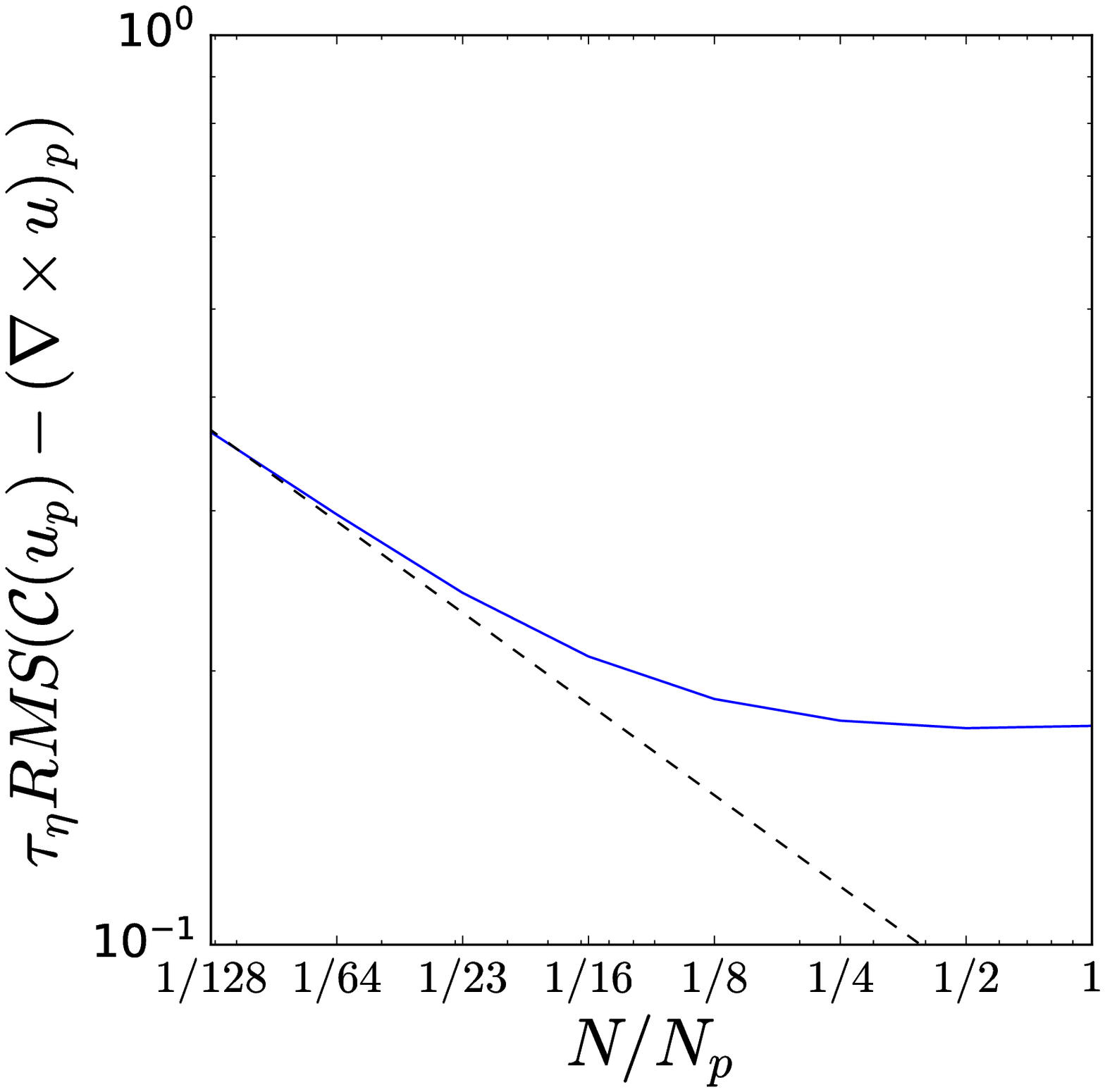}
\caption{PDFs of fluid vorticity normalized by the Kolmogorov time scale $\tau_\eta$ (a) of the fluid vorticity for different values of $N$ computed for fluid particles and on the Eulerian grid. 
RMS error (b) of fluid vorticity computed using the modified Voronoi based method. The dashed line indicates the slope -1/3. 
}
\label{fig:Fluid_Vorticity_VAR}
\end{figure}

\begin{table}
\begin{center}
\begin{tabular}{c c c c c c c c}
 & Fluid & $N_p$ & $N_p/2$ & $N_p/4$ & $N_p/8$ & $N_p/16$ & $N_p/128$  \\
 \hline
$\mathbb{V}(\tau_\eta{\cal C})$ & 0.310 & 0.311 & 0.309 & 0.299 & 0.288 & 0.274 & 0.204  \\
$\mathbb{F}(\tau_\eta{\cal C})$ & 8.831 & 8.960 & 8.875 & 8.746 & 8.570 & 8.334 & 7.268  
\end{tabular}
\end{center}
\caption{Variance $\mathbb{V}$ and flatness $\mathbb{F}$ 
of the curl normalized by the Kolmogorov time scale $\tau_\eta$ of the fluid particles for different particle numbers with $N_p= 1.07 \times 10^9$ and for the fluid vorticity. 
}
\label{tab:var_flat_fluid_curl}
\end{table}

Figure \ref{fig:Fluid_Vorticity_VAR}(a) shows the PDF of the fluid vorticity normalized by the Kolmogorov time scale $\tau_\eta = \nu^{1/2}\epsilon^{-1/2}$ computed at particle positions for different number of particles $N$. Here, $\nu = 1.10\times 10^{-3}$ is the kinematic viscosity and $\epsilon=0.344$ is the energy dissipation rate.
As the 
flow is statistically isotropic, no distinction is made between each component of the vorticity.
We can observe that when the number of particles increases, the variance of the PDF increases. 
For $N \geq N_p/2$, the PDFs of the vorticity computed at the particle position and on the Eulerian grid almost perfectly superimpose. 
Figure \ref{fig:Fluid_Vorticity_VAR}(b) shows the 
root mean square (RMS) error $\sqrt{\sum_{i=0}^{N-1}({\cal C}({ u}_{p_i}) - (\nabla \times u)_{p_i})^2/N}$
of the fluid vorticity computed using the modified Voronoi based method. As shown in 
section~\ref{sec:Computation_of_a_spatial_derivative}, 
we can observe a decrease in the error followed by a plateau. 
From $N \ge N_p/4$ we find that the error saturates and remains constant. We tested if the saturation is not due to the coupling with $\Delta t$, and we conjecture that it is due to linear interpolation of the exact fluid vorticity at the particle position.


%
Table~\ref{tab:var_flat_fluid_curl} assembles the variance and flatness values of the fluid vorticity for different particle numbers. We can observe that the variance and the flatness increase with the number of particles and reaches values close to the DNS values for $N = N_p/2$ and $N = N_p$. 
This number of particle corresponds to an average of one or more particles per Kolmogorov scale, i.e., $N\eta^3/(2\pi)^3 \approx 1$, where $\eta = \nu^{3/4} \epsilon^{-1/4}$ is the Kolmogorov scale. 
%
%

For fluid particles, the Pearson correlation between the vorticity of the DNS data and the vorticity computed at particle positions is $R[{\mathbfcal{C}}(\bm u_p), (\nabla \times \bm u)_p]=0.93$ for $N = N_p/128$ and $R[{\mathbfcal{C}}(\bm u_p), (\nabla \times \bm u)_p]=0.98$ for $N = N_p$.
%
%
%
The correlation obtained for $N = N_p/128$ is 
larger than the correlation coefficient estimated based on the synthetic flow, which is $R[{\mathbfcal{C}}(\bm u_p), (\nabla \times \bm u)_p]=0.46$, as represented by the black square in figure 6(b).
This is because the turbulent flow has 
a faster decay of the spectrum due to dissipation and 
coherent structures, while the synthetic random velocity field tested in figure \ref{fig:pearson_2D_3D_frequence}(b)
has no coherent structures.

\begin{figure}
\centering
\includegraphics[width=0.45\linewidth]{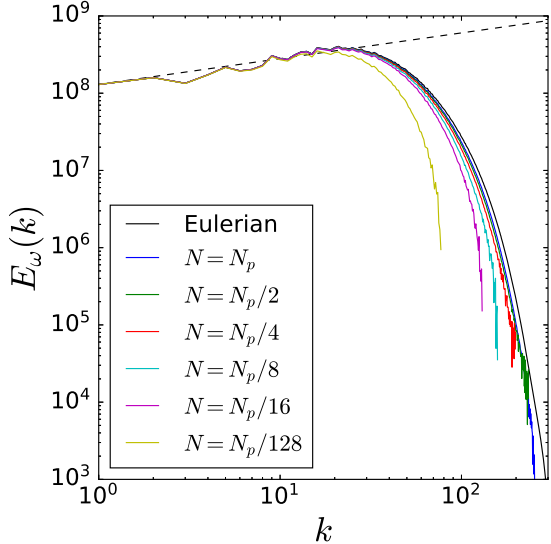}
\caption{Enstrophy spectrum $E_\omega(k)$ of the Eulerian vorticity $\bm \omega$ and of the Lagrangian vorticity $\mathbfcal{C}({\bm u}_p)$ for different particle numbers $N$, all computed on an Eulerian grid with resolution $512^3$. In the spectra of the Lagrangian vorticity the Poisson noise has been removed, by subtracting $C_Nk^2$. The dashed line corresponds to $k^{1/3}$. 
}
\label{fig:Enstrophy_spectrum_fluid_poisson}
\end{figure}

To get insight into the scale distribution of the enstrophy, the squared $L^2$-norm of the vorticity, we show in 
figure \ref{fig:Enstrophy_spectrum_fluid_poisson} the fluid enstrophy spectrum $E_\omega(k)$ 
computed on an Eulerian grid of size $512^3$ and using the modified Voronoi based method for different particle numbers $N$. The dashed line represents the theoretical power law 
of the enstrophy spectrum in HIT, i.e., $k^{1/3}$.
Poisson noise has been removed by subtracting the function $C_Nk^2$ where $C_N$ is a coefficient which depends on the number of particles $N$. 
First, for $N=N_p$ we almost perfectly recover the Eulerian enstrophy spectrum. 
The small deviation at large wavenumbers can be attributed to the effect of the box kernel when the curl value is projected onto the Eulerian grid.
Hence, we can consider that $N=N_p$ is sufficient to capture all scales of motion in the turbulent flow field using the Voronoi based technique.
Reducing the number of particles implies a faster decay of the enstrophy spectra and can thus be interpreted as a low pass filter.
%
The variation in the number of particles can therefore capture the enstrophy at different scales. 


\begin{figure}
\centering
\includegraphics[width=0.45\linewidth]{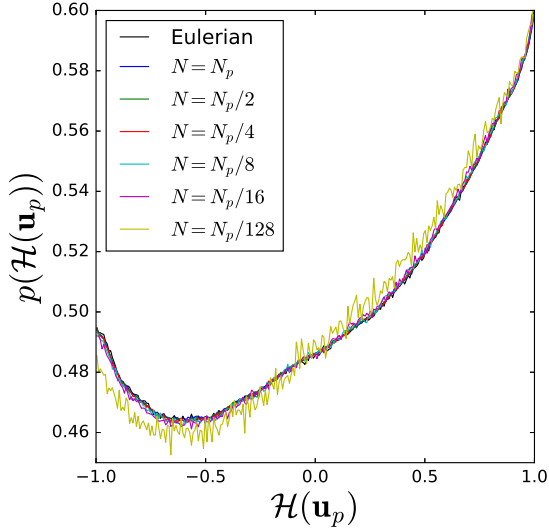}
\caption{PDFs of the relative helicity (cosine of the angle between particle velocity and vorticity) ${\cal H}({\bm u}_p)$ of the particle velocity computed for fluid particles (Lagrangian) and on the Eulerian grid.
}
\label{fig:PDF_Helicity_1G}
\end{figure}

\medskip

Having access to the curl of the particle velocity,
the relative helicity of the particle velocity can be defined as
\begin{equation}
    {\cal H}({\bm u}_p) = \frac{{\bm u}_p\cdot \mathbfcal{C}({\bm u}_p)}{||{\bm u}_p||_2 ~ ||\mathbfcal{C}({\bm u}_p)||_2}, 
\end{equation}
%
which yields the cosine of the angle between the particle velocity and its curl and thus characterizes their alignment. Values of $+1$ correspond to alignment of both quantities, $-1$ to anti-alignment and $0$ to two-dimensional motion, where the curl is perpendicular to the velocity. Helicity is typically used to quantify swirling motion and coherent vortices, see e.g. \citet{moffatt1992helicity}. 
Therefore here we also examine the ability of the proposed method for evaluating the relative helicity. 
Figure \ref{fig:PDF_Helicity_1G} shows the PDF of the relative helicity of fluid particle velocity.
We observe that the helicity PDFs computed on the Eulerian grid and the Lagrangian one computed on the particle positions perfectly superimpose for $N = N_p$.
Although reducing the number of particles has an impact on the amplitude of the vorticity, it has a weak impact 
on the PDF of the relative helicity. It remains almost unchanged
when reducing the number of particles by one order of magnitude, i.e. 
$N \ge N_p/16$, which implies that the alignment between vorticity and velocity remains statistically the same.
However, if we continue to reduce the number of particles, here $N = N_p/128$, we see that the PDF starts to deviate from the exact value.
%


\subsection{Application to inertial particles}
%
%
%
Now we consider one-way coupled inertial heavy point particles \citep{Maxe87} 
for Stokes number of unity in high Reynolds number turbulence.
The Stokes number characterizes the the inertial dynamics of particles and is define as $St = \tau_p / \tau_\eta$, where $\tau_p$ is the particle relaxation time and $\tau_\eta$ the Kolmogorov time.
The equations of the particle position ${\bm x_p}$ and the particle velocity ${\bm v_p}$ are given by
\begin{equation}
    d_t {\bm x_p} = {\bm v_p} \, , \quad \quad d_t {\bm v_p} = - \frac{{\bm v_p} - {\bm u_p}}{\tau_p}  
    \label{eq:particle} \, ,
\end{equation}
where ${\bm u_p}$ is the fluid velocity at particle position ${\bm x_p}$ and time $t$.
For details on the DNS data and the definition of the different parameters we refer to \cite{matsuda2021scale}.
%
%
%
For inertial particles the reference solutions for divergence, curl and helicity of the particle velocity are unknown. 

\begin{figure}
\centering
(a)
\includegraphics[width=0.45\linewidth]{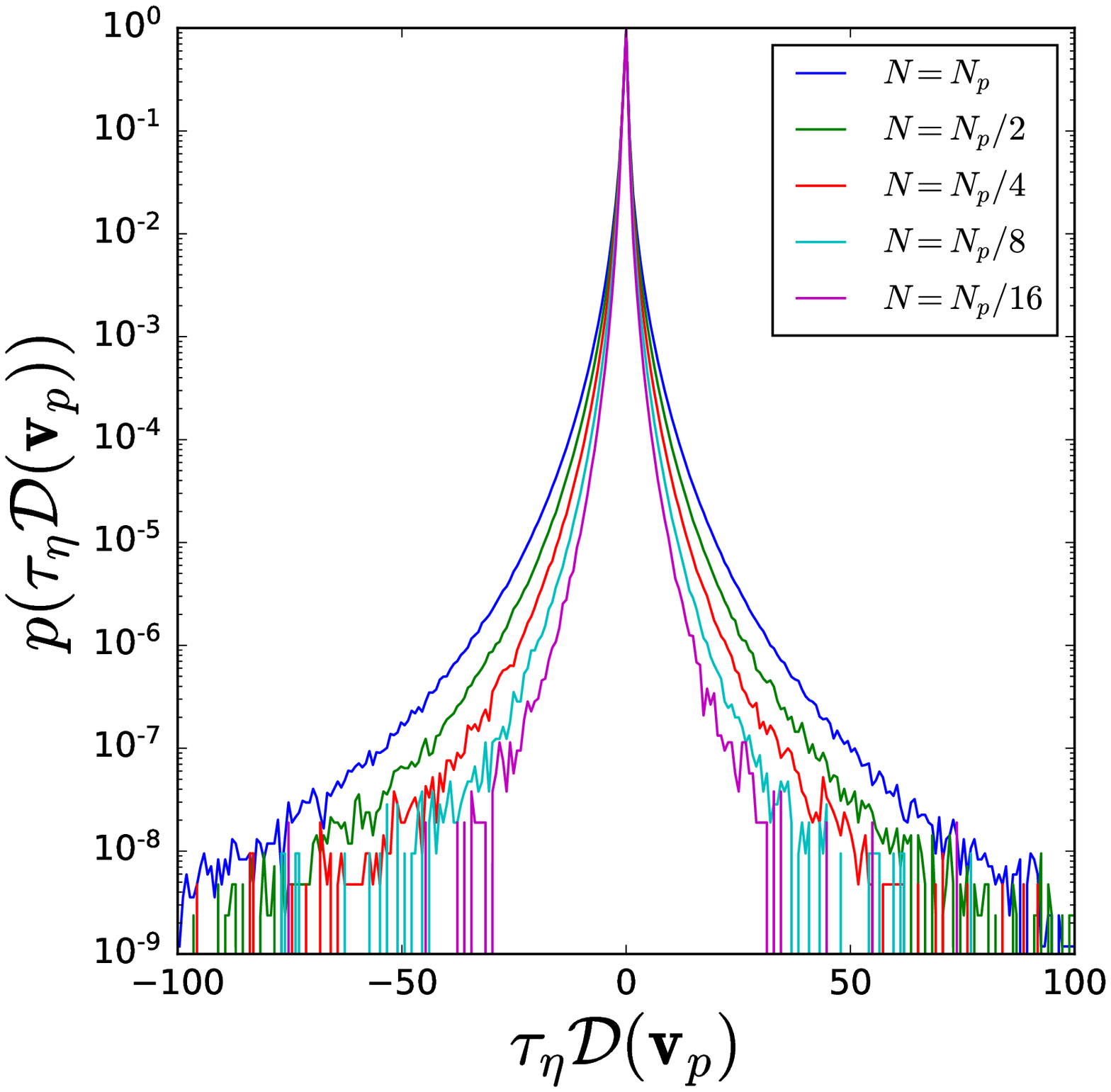}
(b)
\includegraphics[width=0.45\linewidth]{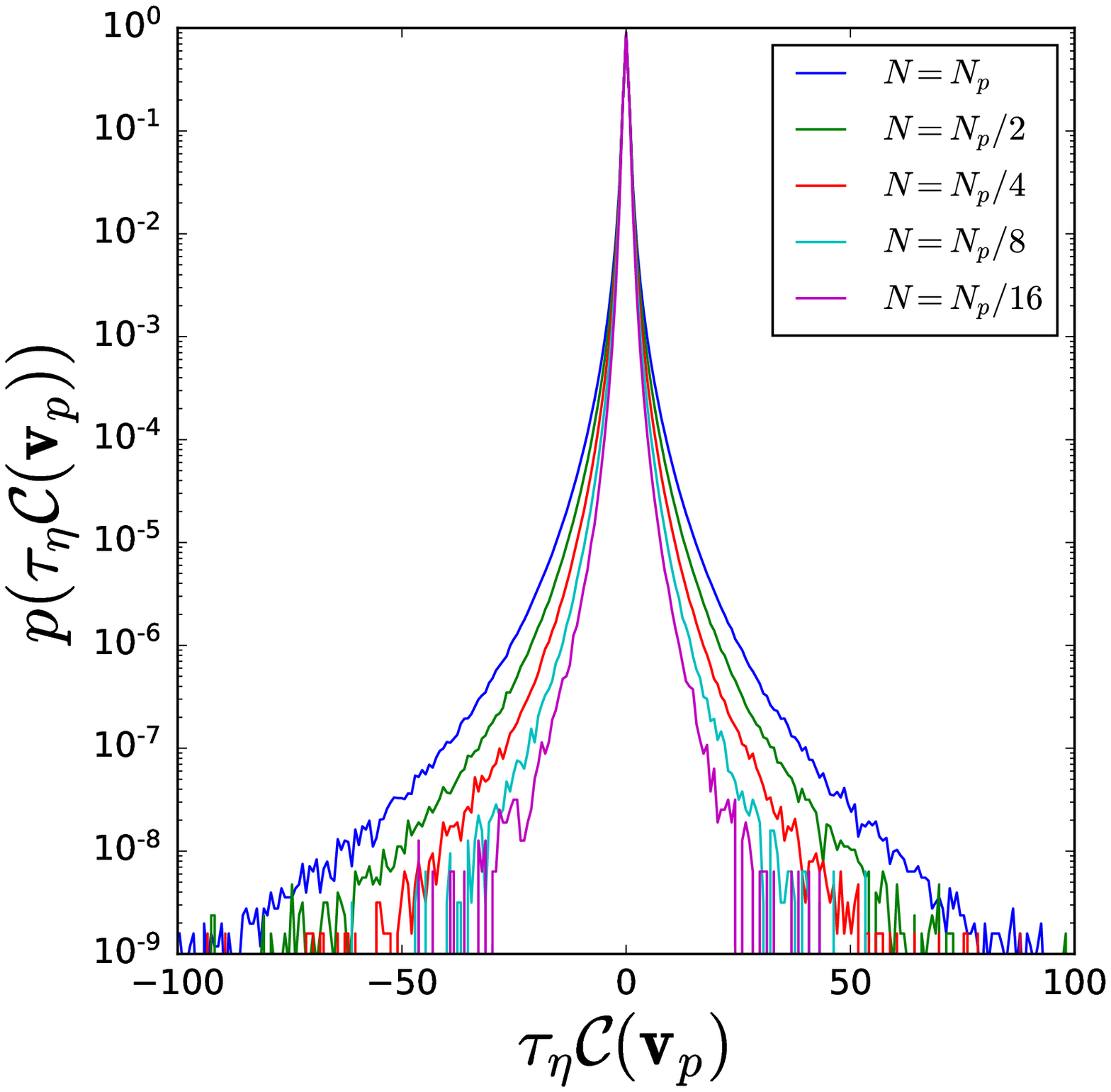}
\caption{PDFs of the divergence ${\cal D}({\bm v}_p)$ (a) and curl $\mathbfcal{C}({\bm v}_p)$ (b)  of the particle velocity computed for inertial particles for $St=1$, normalized by the Kolmogorov time scale $\tau_\eta$. 
}
\label{fig:PDF_Divergence_Curl_St1_1G}
\end{figure}
%
%
%
Figure \ref{fig:PDF_Divergence_Curl_St1_1G} shows the normalized PDF of the divergence and the curl for $St=1$ for different number of particles $N$. 
Reducing the number of particles from 
$N = N_p$
subsequently by a factor of $2$ illustrates that the tails of both divergence and curl PDFs are affected and the range of extreme values is reduced, 
reflected likewise in a change of the variance. 
%
%
%
\begin{table}
\begin{center}
\begin{tabular}{c c c c c c}
 & $N_p$ & $N_p/2$ & $N_p/4$ & $N_p/8$ & $N_p/16$  \\
 \hline
$\mathbb{V}(\tau_\eta{\cal D})$ & 1.321  & 0.896  & 0.610 & 0.418 & 0.287 \\
$\mathbb{V}(\tau_\eta{\cal C})$ & 0.636  & 0.475  & 0.367 & 0.293 & 0.240 \\
$\mathbb{F}(\tau_\eta{\cal D})$ & 126.6  & 117.1  & 130.2 & 83.47 & 92.43 \\
$\mathbb{F}(\tau_\eta{\cal C})$ & 131.6  & 91.32  & 65.76 & 38.45 & 27.13 \\
\end{tabular}
\end{center}
\caption{Variance $\mathbb{V}$ and flatness $\mathbb{F}$ 
of the divergence and the curl normalized by the Kolmogorov time scale $\tau_\eta$ of the particle velocity computed for inertial particles for $St=1$ for different particle numbers with $N_p= 1.07 \times 10^9$. 
}
\label{tab:var_div_curl}
\end{table}
Table~\ref{tab:var_div_curl} assembles the variance values of the divergence and the curl of the particle velocity for different particle numbers. The results show that the variance increases when increasing the number of particles. Fitting a power-law, it follows that
$\mathbb{V}({\cal D}) = C N_p^\alpha$ with $C=1.42\times 10^{-5}$ and $\alpha=5.5\times 10^{-1}$ 
and 
$\mathbb{V}({\cal C}) = C N_p^\alpha$ with $C=4.2\times 10^{-4}$ and $\alpha=3.5\times 10^{-1}$
which indicates that there is no saturation.
The flatness of the divergence oscillates around 100, and for the curl the flatness increases from 27 to 131. This confirms the  strong non Gaussianity for both distributions and their spatial intermittency. 

%
In \citet{matsuda2021scale}, it is observed that inertial particle clustering does not have a smallest scale as implied by the power-law behavior of the number density spectrum. This implies that an infinite number of particles would be necessary for obtaining the convergence. A multiscale analysis, as done in \citet{matsuda_ctr2022}, is therefore required to study the particle dynamics scale by scale.

%

%
%

%
%
%
\begin{figure}
\centering
(a)
\includegraphics[width=0.45\linewidth]{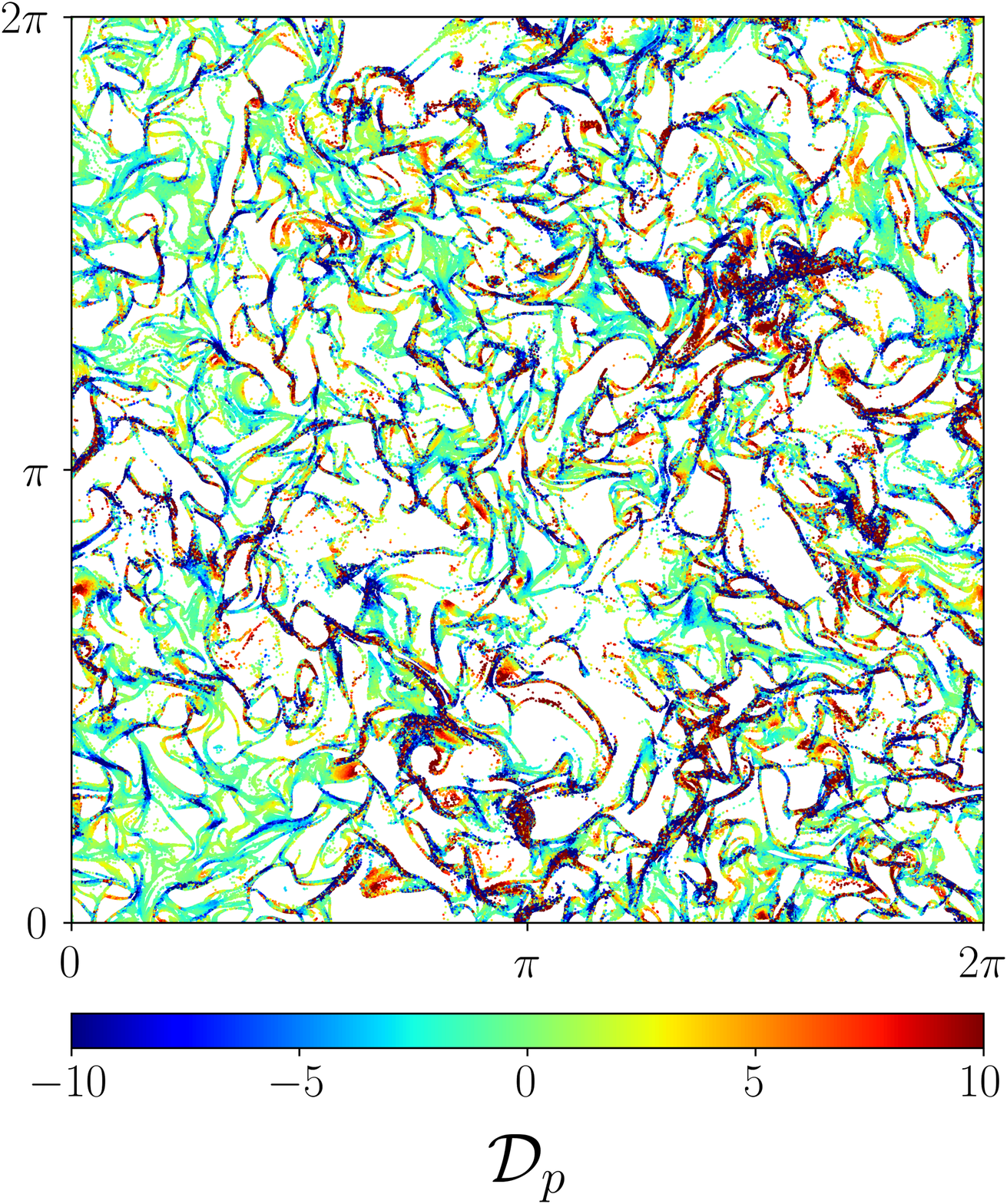}
(b)
\includegraphics[width=0.45\linewidth]{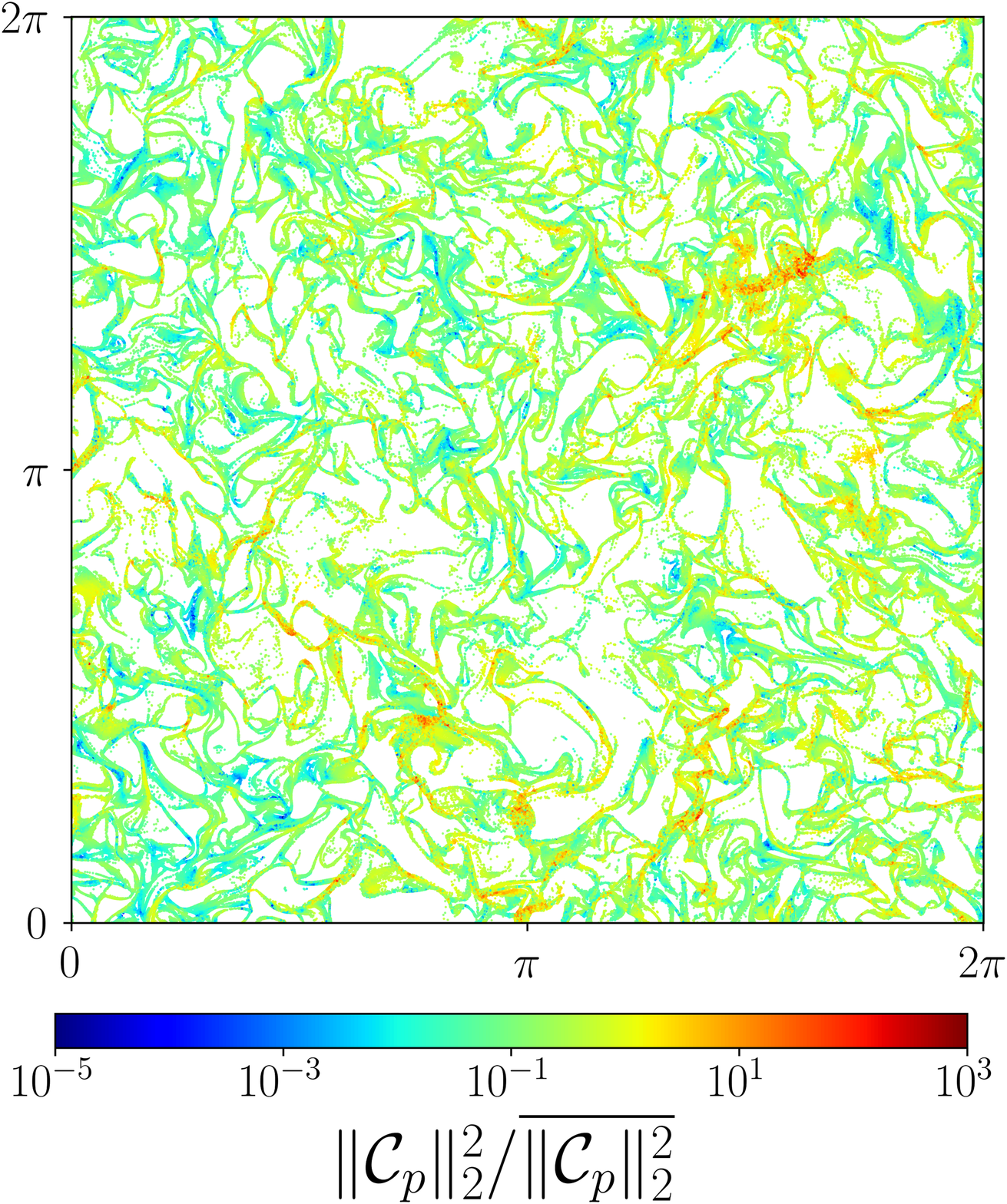}
\caption{Inertial particle positions colored by the divergence (a) and curl (b) for $St=1$ in a slice of $2\pi/512$.
}
\label{fig:Position_Divergence_Curl_1G_St1}
\end{figure}
Note that the finite number of particles is an inherent nature of particle-laden flows and the computed divergence and curl 
values of the inertial particle velocity are only an approximation, as the continuum limit is not known. 
The present method could be used to obtain a low-pass filtered divergence and curl values depending on the number of particles and could give some information about the spatial structures of particle dynamics.
Figure~\ref{fig:Position_Divergence_Curl_1G_St1} shows inertial particle positions colored by particle velocity divergence (a) and curl (b) for $St=1$ in a two-dimensional slice. We can see three different types of clusters in figure~\ref{fig:Position_Divergence_Curl_1G_St1}(a). Clusters with zero divergence where the particles are locally transported, clusters where the sign of the divergence is the same and the transition is continuous, and clusters where the divergence is both positive and negative which is a sign that the particles are crossing. In figure~\ref{fig:Position_Divergence_Curl_1G_St1}(b), we observe that the extreme values of the curl are in the same region where the divergence is both positive and negative which corresponds to strong rotation of the particle clusters.

\begin{figure}
\centering
\includegraphics[width=0.45\linewidth]{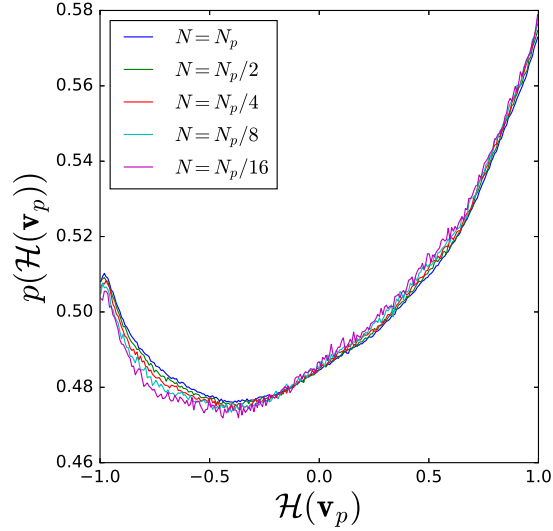}
\caption{PDFs of the relative helicity ${\cal H}({\bm v}_p)$ of the particle velocity computed for $St=1$.
}
\label{fig:PDF_Helicity_1G_St1}
\end{figure}
To quantify the swirling motion of inertial particles, we consider the helicity of the particle velocity.
The PDF of the cosine of particle velocity and its curl is plotted in figure~\ref{fig:PDF_Helicity_1G_St1}. 
The behavior is qualitatively similar to what has been observed for fluid particles in figure~\ref{fig:PDF_Helicity_1G}. Similarly we find that reducing the number of particles 
only weakly impacts the PDFs.

\begin{figure}
\centering
\includegraphics[width=0.45\linewidth]{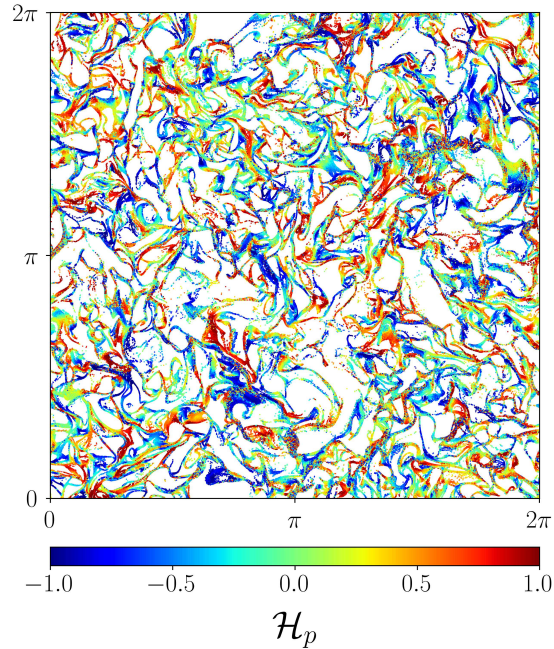}
\caption{Particle positions colored by the helicity for $St=1$ in a slice of $2\pi/512$.
}
\label{fig:Position_Helicity_1G_St1}
\end{figure}
Visualizations of the inertial particles colored with the relative helicity value, given in figure~\ref{fig:Position_Helicity_1G_St1}, allows us to detect and quantify that some clusters exhibit strongly helical motion, while some are not swirling and satisfy quasi planar motion.
We can observe that the value of helicity within a cluster is generally of the same sign, i.e., within a cluster there is no great variation in the value of helicity 
contrary to what is seen for the divergence.

\section{Conclusions}
\label{sec:conclusions}

We proposed tessellation-based methods to determine first order differential operators of the particle velocity advecting particle clouds based on the method proposed in \citet{oujia_matsuda_schneider_2020}. 
The underlying idea is to compute the time change of the volume change rate of the Voronoi volume obtained by applying a tessellation of the particle positions.
We proposed a modification of the Voronoi tessellation, i.e., by using the centers of gravity of the Delaunay cells as the vertices of the modified Voronoi cells. We showed that the modified method converges with first order in time and in each space direction for randomly distributed particles, which is an improvement compared to the method proposed in \citet{oujia_matsuda_schneider_2020}. 
%
%
%
%
We showed that the method originally developed to compute the divergence of the particle velocity can be likewise used to compute the curl and velocity gradient tensor by rearranging the velocity vector coefficients.
More generally, 
we notice that the method has the ability to compute simultaneously the sum of several first order derivatives. 
%
The modified Voronoi-based method is also applied to synthetic turbulent random flow fields with imposed power-law energy spectra. 
A higher correlation of the Pearson coefficient than what is obtained at the highest wavenumber is shown 
which can be attributed 
to the rapid decay ($k^{-3}$ in two dimensions and $k^{-5/3}$ in three dimensions) of the energy spectrum with wavenumber.
%

Then, we applied the method to fluid particles advected in three-dimensional fully developed turbulence and showed that the method gives accurate results for the computation of vorticity and helicity of the fluid particle velocity.
A strong Pearson correlation is observe for fluid vorticity between the exact and computed value. It is also observed that for a given number of particles, the Pearson correlation is higher for DNS data than for synthetic turbulent random flow fields, which can be ascribed to the fact that the synthetic turbulence as no coherent structures and that DNS has a faster decay of the spectrum due to dissipation. 
The enstrophy spectrum for different particle numbers reveals that a decrease in particle number results in a faster decay of the enstrophy spectra, acting thus as a low pass filter. This implies that varying the number of particles can capture enstrophy at different scales, 
which paves the way for multiscale analysis
using the modified Voronoi-based method.
Preliminary results on the development of a multiscale modified Voronoi tessellation can be found in \citet{matsuda_ctr2022}.
Eventually, for inertial particles in turbulence we showed that the introduced methods for computing the divergence/convergence, the curl and the helicity of the particle velocity allow to assess the clustering dynamics and the vortical and swirling motion of particle clouds, which is of great benefit for analyzing turbulent multiphase flow.

For instance, in \citet{west_ctr2022}, we applied the above technique to four-way coupled turbulent channel flow data, laden with inertial particles and analyzed the influence of mass loading and inertia of the particles on their clustering in the different flow layers. We computed the divergence and the curl in the log-layer, the buffer layer and the viscous sub-layer and quantified the flow behavior.
we studied particle laden isotropic turbulence investigating divergence and curl of the inertial particle velocity. We found extreme divergence and rotation values reflecting the strong intermittency of the particle motion.

Finally, let us mention that
this method based on the variation of the particle position at two time instants seems to be a natural approach to study the particle motion of PTV data. In addition, the method to compute the divergence from the time variation of the volume can be used in fields where cells appear naturally as in biology, or in astronomy for quantifying the clustering of stars in interstellar clouds.
Moreover, the proposed methods can be applied to high dimensional data, e.g. point clouds in six dimensions encountered in plasma physics governed by kinetic equations, or in data science in general.

\section*{Acknowledgments}
\label{sec:ACK}
T.O. and K.S. acknowledge funding from the Agence Nationale de la Recherche (ANR), grant ANR-20-CE46-0010-01.
K.M acknowledges financial support from JSPS KAKENHI Grant Number JP20K04298. 
Centre de Calcul Intensif d’Aix-Marseille is acknowledged for granting access to its high performance computing resources.
The DNS data analyzed in this project were obtained using the Earth Simulator supercomputer system of JAMSTEC.
The authors acknowledge use of computational resources from the Yellowstone cluster awarded by the National Science Foundation to the Center for Turbulence Research (CTR) at Stanford.

\hfill

\newpage



\end{document}